\documentclass[12pt]{article}
\usepackage{amsfonts}
\font\ab=msbm10 scaled\magstep1%Extra symbols, gr 2, medium weight (mathbb)
\newcommand{\R}{\mbox{\ab R}}%R double
\newcommand{\N}{\mbox{\ab N}}%Z double
\newtheorem{thm}{Theorem}[section]
\newtheorem{rem}[thm]{Remark}

\newtheorem{cor}[thm]{Corollary}
\newtheorem{prop}[thm]{Proposition}

\catcode `\@=11
\@addtoreset{equation}{section}

\def\harr#1#2{\smash{\mathop{\hbox to .5in{\rightarrowfill}}
\limits^{\scriptstyle#1}_{\scriptstyle#2}}}

\def\harrl#1#2{\smash{\mathop{\hbox to .5in{\leftarrowfill}}
\limits^{\scriptstyle#1}_{\scriptstyle#2}}}

\def\qed{\vrule height 6pt width 6pt depth 6pt}

\newcommand{\be}{\begin{equation}}
\newcommand{\ee}{\end{equation}}

\begin{document}
\begin{titlepage}
%\thispagestyle{empty}
%\begin{flushright}
%IFA-FT-402-1994, November
%\end{flushright}
%\bigskip\bigskip
\begin{center}
{\bf \Large{On a Generalisation of the Poincar\'e-Cartan Form to Classical 
Field Theory}}
\end{center}
\vskip 1.0truecm
\centerline{D. R. Grigore
\footnote{e-mail: grigore@theor1.ifa.ro, grigore@roifa.ifa.ro}}
\vskip5mm
\centerline{Dept. of Theor. Phys., Inst. Atomic Phys.}
\centerline{Bucharest-M\u agurele, P. O. Box MG 6, ROM\^ANIA}
\vskip 2cm
\bigskip \nopagebreak
\begin{abstract}
\noindent
We present here a possible generalisation of the Poincar\'e-Cartan form in
classical field theory in the most general case: arbitrary dimension,
arbitrary order of the theory and in the absence of a fibre bundle structure.
We use for the kinematical description of the system the $(r,n)$-Grassmann 
manifold associated to a given manifold $X$, i.e. the manifold of $r$-contact
elements of $n$-dimensional submanifolds of $X$. The idea is to define globally
a $n+1$ form on this Grassmann manifold, more precisely its class with respect
to a certain subspace and to write it locally as the exterior derivative of a
$n$ form which is the Poincar\'e-Cartan form.  As an important application we
obtain a new proof for the most general expression of a variationally trivial
Lagrangian.
\end{abstract}



\newpage
\setcounter{page}1
\end{titlepage}

\section{Introduction}

It is widely accepted that the variational principles should be given in a
coordinate independent formulation. This idea was first realised for a
dynamical system with a finite number of degrees of freedom (i.e. particle
mechanics), using a differential 1-form instead of the Lagrangian, by
Poincar\'e and Cartan \cite{Po}, \cite{Ca}. There are a number of
generalisations of this idea for classical field theory \cite{Kr1}, 
\cite{Be1}, \cite{Be2}, \cite{Ru1}, \cite{Ru2}, \cite{Ki}, \cite{GS}, 
\cite{Ga}, \cite{Go}, \cite{Sa}, \cite{Ol2}. A related concept is that of
Lepage equivalent of a Lagrangian form (see for instance \cite{Kr2},
\cite{Kr3}, \cite{Kr9}). All these generalisations use as geometric framework
for classical field theory the jet bundle formalism (more explicitly the
space-time and field variables are local coordinates on a fibre bundle $X$ over
a ``space-time manifold" $M$) and the derivative of the fields, up to order
$s$, are variables in the $s$-th order jet bundle extension $J^{s}X$ of $X$.

It was later \cite{Kl}, \cite{Sou} suggested that it is more convenient to work
with the exterior differential of the above Poincar\'e-Cartan form. For the
case of finite number of degrees of freedom this 2-form is in general 
presymplectic and was used by Souriau and others \cite{Sou}, \cite{Horv} to 
obtain the phase space as a symplectic manifold in in a deductive way.
The idea is to consider that the fundamental  mathematical object for a
Lagrangian system must be this 2-form and not the Lagrangian function or the
Poincar\'e-Cartan 1-form. This point of view leads to the main features of the
Lagrangian and the Hamiltonian formalism and to a natural definition of the 
Noetherian symmetries. 

One can generalise this Lagrange-Souriau form without using the fibration 
hypothesis mentioned above in two particular but important cases: for classical
field theory of first order \cite{Gr1} and for systems with a finite number of
degrees of freedom and of arbitrary order \cite{Gr2}. Moreover, this
Lagrange-Souriau form can be locally written as the exterior differential of a
Poincar\'e-Cartan form related to some chosen local chart. This
Poincar\'e-Cartan form is the same as that given by Krupka \cite{Kr1}, Betounes
\cite{Be1}-\cite{Be2} and Rund \cite{Ru1}. The Lagrangian is locally determined
up to a variationally trivial Lagrangian, i.e. a Lagrangian giving trivial
Euler-Lagrange equations. As a consequence, one can define in a geometrically
nice way the Noetherian symmetries using the Lagrange-Souriau form. 

In this paper we give a generalization of the Lagrange-Souriau and of the
Poincar\'e-Cartan forms in the most general case used in classical field
theory. We consider an arbitrary manifold $X$ without a fibre bundle structure
over some space-time manifold so instead of the $s$-th order jet bundle
extension one must use the $s$-th order Grassmann bundle $P^{s}_{n}X$ 
associated to $X$ which was recently considered in the literature \cite{GK}.
In the next Section we will summarise the main features of this construction.
Next, in Section 3, we will be able to define globally a $n+1$ differential
form, but one will be able to see that, in general, one cannot determine this
form uniquely.  Fortunately one can consider the equivalence class of this form
to a certain globally defined subspace of differential forms. This equivalence
class is the ``physical" object we are looking for. It is interesting to note
that this subspace of differential forms is in fact identically zero exactly 
in the two particular cases mentioned above ($s = 2$, $n$ arbitrary and $n =
1$, $s$ arbitrary). We will need some combinatorial tricks introduced in
\cite{Gr3} to simplify the analysis of some tensorial identities. In Section 4
we locally exhibit the $n+1$ differential form as the exterior derivative of
a locally defined Poincar\'e-Cartan $n$-form and in this way the (local)
Lagrangian function appears also. We also present in this Section some natural
set of equations having solutions in terms of the so-called hyper-Jacobians
\cite{BCO}, \cite{Ol1}. In fact, we consider that this set of equations, which
emerges naturally from a closedness condition imposed on some differential
form, and its natural solution in terms of hyper-Jacobians is a key point of
this paper. Next, we use the formalism developed above to provide the
most general expression for a variationally trivial Lagrangian of arbitrary
order already obtained in \cite{Gr4} by a different method. Finally, in Section
5 we present two particular but very important cases, namely $n = 1$ $s$
arbitrary and $s = 2$ $n$ arbitrary. Some ideas related to the ones from this
paper also appear in \cite{Kr5}.

\newpage

\section{Grassmann Manifolds}

\subsection{The Basic Constructions of the Grassmann Manifolds}

In this Section we present the basic construction of Grassmannian manifolds
following \cite{GK} and \cite{Gr5}. We will skip all the proofs.
We consider
$N$, $n\ge 1$
and
$r\ge 0$
integers such that
$n\le N$,
and let
$X$
be a smooth manifold of dimension
$N$ which is the mathematical model for the kinematical degrees of freedom of a
certain classical field theory. 

Let
$U \subset \R^{n}$
be a neighbourhood of the point 
$0 \in \R^{n}$,
$x \in X$
and let
$\Gamma_{(0,x)}$
be the set of smooth immersions 
$\gamma: U \rightarrow X$
such that 
$\gamma(0) = x$.
On 
$\Gamma_{(0,x)}$
one has the the equivalence relationship
$``\gamma \sim \delta"$
{\it iff} there exists a chart 
$(V,\psi) \quad \psi = (x^{A}), \quad A =1, \dots ,N$
on $X$ such that the functions
$\psi \circ \gamma, \psi \circ \delta: \R^{n} \rightarrow \R^{N}$
have the same partial derivatives up to order $r$ in the point $0$.
The equivalence class of
$\gamma$
will be denoted by
$j^{r}_{0}\gamma$
and it is called a $(r,n)$-{\it velocity}. The set of
$(r,n)$-velocities
at
$x$
is denoted by
$T^{r}_{(0,x)}(\R^{n},Y) \equiv \Gamma_{(0,x)} /\sim$.
We denote
$$
T^{r}_{n}X = \bigcup_{x \in X} T^{r}_{(0,x)}(\R^{n},X), 
$$
and define surjective mappings
$\tau^{r,s}_{n}: T^{r}_{n}X \rightarrow T^{s}_{n}X$, 
where
$0 < s \leq r$,
by
$
\tau^{r,s}_{n}(j^{r}_{0}\gamma) = j^{s}_{0}\gamma
$
and
$\tau^{r,0}_{n}: T^{r}_{n}X \rightarrow X$, 
where
$1 \leq r$,
by
$
\tau^{r,0}_{n}(j^{r}_{0}\gamma) = \gamma(0).
$

If
$(V,\psi), \quad \psi =(x^{A})$,
is a chart on
$X$
we define the couple 
$(V^{r}_{n},\psi^{r}_{n})$ 
where 
$V^{r}_{n} = (\pi^{r,0}_{n})^{-1}(V),$
$
\psi^{r}_{n} = (x^{A},x^{A}_{j}, \cdots, x^{A}_{j_{1},j_{2},\dots,j_{r}})
\quad 1\leq j_{1} \leq j_{2} \leq \cdots \leq j_{r} \leq n,
$
and
\be
x^{A}_{j_{1},\dots,j_{k}}(j^{r}_{0}\gamma) \equiv 
\left. {\partial^{k} \over \partial t^{j_{1}} \dots \partial t^{j_{k}}}
x^{A} \circ \gamma\right|_{0}, 
\quad 0 \le k\le r.
\label{coord-J}
\ee

The expressions
$x^{A}_{j_{1}, \cdots j_{k}}(j^{r}_{0}\gamma)$
are defined for all indices
$j_{1},\dots ,j_{r}$ in the set 
$\{1,\dots ,n\}$
but because of the symmetry property
\be
x^{A}_{j_{P(1)},\dots,j_{P(k)}}(j^{r}_{0}\gamma) =
x^{A}_{j_{1},\dots,j_{k}}(j^{r}_{0}\gamma) \quad (k = 2,...,n)
\ee
for all permutations
$P \in {\cal P}_{k}$
of the numbers
$1, \dots ,k$
we consider only the independent components given by the restrictions
$
1\leq j_{1} \leq j_{2} \leq \cdots \leq j_{r} \leq n.
$
This allows one to use multi-index notations i.e.
$
\psi^{r}_{n} = (x^{A}_{J}), \quad |J| = 0,...,r
$
where by definition
$x^{A}_{\emptyset} \equiv x^{A}.$
The same comment is true for the partial derivatives
${\partial \over \partial x^{A}_{j_{1},\dots,j_{k}}}$.

The couple 
$(V^{r}_{n},\psi^{r}_{n})$
is a chart on 
$T^{r}_{n}X$
called the {\it associated chart} of the chart
$(V,\psi)$
and the system of charts give a smooth structure on this set;
moreover 
$T^{r}_{n}X$
is a fibre bundle over $X$ with the canonical projection
$\tau^{r,0}_{n}$.
The set
$T^{r}_{n}Y$
endowed with the smooth structure defined by the associated charts defined
above is called the {\it manifold of}
$(r,n)$-{\it velocities}
over
$X$.

In the chart
$(V^{r}_{n},\psi^{r}_{n})$
one introduces the following differential operators: 

\be
\Delta^{j_{1},\dots,j_{k}}_{A} \equiv {r_{1}! \dots r_{n}! \over k!}
{\partial \over \partial x^{A}_{j_{1},\dots,j_{k}}} , 
\qquad j_{1},\dots,j_{k} \in \{1, \dots, n\}
\label{partial}
\ee
where 
$r_{k}$
is the number of times the index $k$ shows up in the sequence
$j_{1}, \dots j_{k}$.

The combinatorial factors are such that the following relation is true:

\be
\Delta^{i_{1},\dots,i_{k}}_{A} x^{B}_{j_{1},\dots,j_{l}} =
\cases{ \delta^{B}_{A} {\cal S}^{+}_{j_{1},\dots,j_{k}} \delta^{i_{1}}_{j_{1}}
\dots \delta^{i_{k}}_{j_{k}} & if $k = l$ \cr 0 & if $k \not= l$ \cr}.
\label{derivation}
\ee

Here we use the notations from \cite{Gr2}, namely
${\cal S}^{\pm}_{j_{1},\dots,j_{k}}$
are the symmetrization (for the sign $+$) and respectively the 
antisymmetrization (for the sign $-$) projector operators defined by

\be
{\cal S}^{\pm}_{j_{1},\dots,j_{k}} f_{j_{1},\dots,j_{k}} \equiv
{1 \over k!} \sum_{P \in {\cal P}_{k}} \epsilon_{\pm}(P) 
f_{j_{P(1)},\dots,j_{P(k)}}
\label{sa}
\ee
where the sum runs over the permutation group 
${\cal P}_{k}$
of the numbers
$1, \dots, k$
and
$$
\epsilon_{+}(P) \equiv 1, \quad \epsilon_{-}(P) \equiv (-1)^{|P|}, \quad
\forall P \in {\cal P}_{k};
$$
here
$|P|$
is the signature of the permutation $P$.

In this way one takes care of overcounting the indices. More precisely, for any
smooth function on $V^{r}$, the following formula is true:

\be
df = \sum_{k=0}^{r} (\Delta^{j_{1},\dots,j_{k}}_{A} f) 
dx^{A}_{j_{1},\dots,j_{k}} = 
\sum_{|I| \leq r} (\Delta^{I}_{A} f) dx^{A}_{I}
\label{df}
\ee
where we have also used the convenient multi-index notation.

The {\it formal derivatives} are:

\be
D^{r}_{i} \equiv \sum_{k=0}^{r-1} x^{A}_{i,j_{1},\dots,j_{k}} 
\Delta^{j_{1},\dots,j_{k}}_{A} = 
\sum_{|J| \leq r-1} x^{A}_{iJ} \Delta^{J}_{A}.
\label{formal}
\ee

The last expression uses the multi-index notation; if $I$ and $J$ are two such
multi-indices we mean by $IJ$ the juxtaposition of the two sets $I, J.$

When no danger of confusion exists we simplify the notation putting simply
$D_{i} = D^{r}_{i}.$

The formal derivatives give a conveniently expression for the change of charts
on the velocity manifold induced by a change of charts on $X$. Let
$(V,\psi)$
and
$(\bar{V},\bar{\psi})$
two charts on $X$ such that 
$V \cap \bar{V} \not= \emptyset$
and let
$(V^{r},\psi^{r})$
and
$(\bar{V}^{r},\bar{\psi}^{r})$
the corresponding attached charts from
$T^{r}_{n}X.$
The change of charts on $X$ is
$F: \R^{N} \rightarrow \R^{N}$
given by:
$F \equiv \bar{\psi} \circ \psi^{-1}.$
It is convenient to denote by 
$F^{A}: \R^{N} \rightarrow \R$
the components of $F$ given by
$F^{A} \equiv \bar{x}^{A} \circ \psi^{-1}$.
We now consider the change of charts on 
$T^{r}_{n}X$
given by
$F^{r} \equiv \bar{\psi}^{r} \circ (\psi^{r})^{-1}.$
One notes that
$V^{r} \cap \bar{V}^{r} \not= \emptyset$;
we need the explicit formulae for the components of 
$F^{r}$, namely for the functions
$$
F^{A}_{j_{1}, \dots, j_{k}} \equiv \bar{x}^{A}_{j_{1}, \dots, j_{k}} \circ
(\psi^{r})^{-1}, \quad j_{1} \leq j_{2} \dots \leq j_{k}, \quad k = 1,...,r
$$
defined on the overlap:
$V^{r} \cap \overline V^{r}$.

One finds out that the functions
$F^{A}_{j_{1}, \dots, j_{k}}$
are given by the following relations:

\be
F^{A}_{jI} = D_{j} F^{A}_{I} \quad |I| \leq r-1;
\label{F_I}
\ee
or more explicitly:

\be
F_{I}^{A} = \sum_{p=1}^{|I|} 
{\sum\limits_{(I_{1},\ldots,I_{p})}} 
x_{I_{1}}^{B_{1}} \cdots x_{I_{p}}^{B_{p}}
(\Delta_{B_{1}} \cdots \Delta_{B_{p}} F^{A}), \quad 1 \leq |I| \leq r
\label{change-chart}
\ee
where the second sum denotes summation over all partitions
${\cal P}(I)$
of the set
$I$
and two partitions are considered identical if they differ only by a
permutation of the subsets. 

By definition the {\it differential group of} order $r$ is the set

\be
L^{r}_{n} \equiv \{ j^{r}_{0} \alpha \in J^{r}_{0,0}(\R^{n},\R^{n}) \vert
\alpha \in Diff(\R^{n}) \}
\ee
i.e. the group of invertible
$r$-jets
with source and target at
$0\in \R^{n}$.
The group multiplication in
$L^{r}_{n}$
is defined by the jet composition
$
L^{r}_{n}\times L^{r}_{n} \ni (j^{r}_{0}\alpha, j^{r}_{0}\beta) 
\mapsto  j^{r}_{0}(\alpha\circ \beta) \in L^{r}_{n}.
$

The {\it canonical} (global) {\it coordinates} on 
$L^{r}_{n}$
are defined by

\be
a^{i}_{j_{1},\dots,j_{k}}(j^{r}_{0}\alpha) = \left.
{\partial^{k} \alpha^{i} \over 
\partial t^{j_{1}} \dots \partial t^{j_{k}}}\right|_{0}, \quad
j_{1} \leq j_{2} \leq \dots \leq j_{k},  \quad k = 0,...,r
\label{coord-L}
\ee
where
$\alpha^{i}$
are the components of a representative
$\alpha$
of
$j^{r}_{0}\alpha$.

We denote
$$
a \equiv (a^{i}_{j},a^{i}_{j_{1},j_{2}}, \dots ,a^{i}_{j_{1},\dots,j_{k}}) = 
(a^{i}_{J})_{1 \leq |J| \leq r}
$$
and notice that one has

\be
det(a^{i}_{j}) \not= 0.
\label{group}
\ee 

The composition law for the differential group can be obtained explicitly:

\be
(a \cdot b)^{k}_{I} = \sum^{|I|}_{p=1} {\sum\limits_{(I_{1},\dots,I_{p})}}
b^{j_{1}}_{I_{1}} \dots b^{j_{p}}_{I_{p}} a^{k}_{j_{1}, \dots j_{p}}, 
\quad |I| = 1, \dots r. 
\label{composition}
\ee

The group
$L^{r}_{n}$
is a Lie group.

The manifolds of
$(r,n)$-velocities
$T^{r}_{n}Y$
admits a (natural) smooth right action of the differential group
$L^{r}_{n}$,
defined by the jet composition

\be
(x\cdot a)^{A}_{I} \equiv x^{A}_{I}(j_{0}^{r}(\gamma \circ \alpha)) 
\label{x-a}
\ee
where the connection between 
$x^{A}_{I}$
and
$\gamma$
is given by (\ref{coord-J}) and the connection between 
$a^{i}_{I}$
and
$\alpha$
is given by (\ref{coord-L}). 

The chart expression of this action can also be obtained explicitly:

\be
(x\cdot a)^{A} = x^{A},\quad 
(x\cdot a)_{I}^{A} = 
\sum_{p=1}^{|I|} {\sum\limits_{(I_{1},\dots,I_{p}) \in {\cal P(I)}}} 
a_{I_{1}}^{j_{1}} \dots a_{I_{p}}^{j_{p}} x_{j_{1},\dots,j_{p}}^{A},
\quad |I| \geq 1
\label{action}
\ee
and it is smooth.

The group
$L^{r}_{n}$
has a natural smooth left action on the set of smooth real functions defined on
$T^{r}_{n}X$ ,
namely for any such function $f$ we have:

\be
(a\cdot f)(x) \equiv f(x\cdot a).
\label{act}
\ee

We say that a
$(r,n)$-velocity
$j_{0}^{r}\gamma \in T_{n}^{r}X$
is {\it regular}, if 
$\gamma$
(or any other representative) is an immersion. If
$(V,\psi)$,
$\psi =(x^{A})$,
is a chart, and the target
$\gamma(0)$
of an element
$j_{0}^{r}\gamma \in T_{n}^{r}X$
belongs to
$V$,
then
$j_{0}^{r}\gamma$
is regular {\it iff} there exists a subsequence 
$ {\bf I} \equiv (i_{1},\dots,i_{n})$ 
of the sequence
$(1,2,\dots,n,n+1,\dots,n+m)$
such that

\be
{\rm det} (x^{i_{k}}_{j}) \not= 0;
\label{regular}
\ee
(here 
$x^{i_{k}}_{j}$
is a 
$n \times n$
real matrix.) The associated charts have the form
$$
(V^{{\bf I},r},\psi^{{\bf I},r}), \quad
\psi^{{\bf I},r} = (x^{k}_{I},x^{\sigma}_{I}), \quad k = 1,\dots,n, 
\quad \sigma = 1,\dots m \equiv N - n, \quad |I| \leq r
$$
where
$$
x^{k}_{I} \equiv x^{i_{k}}_{I}, \quad k = 1,\dots n
$$ 
and
$\sigma \in \{ 1,\dots, N\} - \{i_{1},\dots,i_{n}\}.$
The set of regular
$(r,n)$-velocities
is an open,
$L_{n}^{r}$-invariant
subset of
$T_{n}^{r}X$,
which is called the {\it manifold of regular} $(r,n)$-{\it velocities},
and is denoted by
${\rm Imm} T_{n}^{r}X$.

One can find out a complete system of
$L_{n}^{r}$-invariants
of the action (\ref{action}) on
${\rm Imm} T_{n}^{r}X$;
for simplicity we take the chart for which one has
$\{i_{1},\dots,i_{n}\} = \{ 1,\dots, n\}$
and we will denote
$$
x^{\sigma}_{I} \equiv x^{n+\sigma}_{I}, \quad 
\sigma = 1,\dots m, \quad |I| \leq r.
$$

Let
$(V,\psi), \psi = (x^{A})$,
be a chart on
$X$
and let
$(V_{n}^{r},\psi_{n}^{r})$
be the associated chart on
${\rm Imm} T_{n}^{r}X$
with coordinates:
$(x^{\sigma}_{I},x^{i}_{I}).$
Because
${\bf x} \equiv (x^{i}_{I})$
is an element of
$L_{n}^{r}$
we can define its inverse
${\bf z} \equiv (z^{i}_{I})$.
Now we define recurringly on this chart the following functions

\be
y^{\sigma} \equiv x^{\sigma}, \quad
y_{i_{1},\dots, i_{k}}^{\sigma} = z_{i_{1}}^{j} D_{j}
y_{i_{2},\dots, i_{k}}^{\sigma},\quad k = 1,\dots, r;
\label{invariants}
\ee
(here 
$z^{j}_{i}$
are the first entries of the element 
${\bf z} \equiv {\bf x}^{-1}$ 
from
$L^{r}_{n}$.)

Then the functions
$y^{\sigma}_{i_{1},\dots,i_{k}}$
so defined have a number of important properties:
\begin{itemize}
\item
They depend smoothly only on
$x^{A}_{J}, \quad |J| \leq k$;
\item
Are completely symmetric in all indices
$i_{1},\dots, i_{k}, \quad k = 1,\dots r$;
\item
Are {\bf uniquely} determined by the recurrence relations:

\be
x_{I}^{\sigma} = \sum_{p=1}^{|I|} \sum\limits_{(I_{1},\dots,I_{p})}
x_{I_{1}}^{j_{1}} \dots x_{I_{p}}^{j_{p}} y_{j_{1},\dots,j_{p}}^{\sigma}
= (y\cdot {\bf x})^{\sigma}_{I}, \quad 1 \leq |I| \leq r.
\label{yx}
\ee
\item
Because one can ``invert" these formulae to

\be
y^{\sigma}_{I} = (x\cdot {\bf z})^{\sigma}_{I} 
= \sum_{p=1}^{|I|} \sum\limits_{(I_{1},\dots,I_{p})}
z_{I_{1}}^{j_{1}} \dots z_{I_{p}}^{j_{p}} x_{j_{1},\dots,j_{p}}^{\sigma},
\quad 1 \leq |I| \leq r
\label{xy}
\ee
one can use on
$V^{r}$
the new coordinates
$(y^{\sigma}_{I},x^{i}_{I}), \quad |I| \leq r$.
\item
The functions
$y^{\sigma}_{I}, \quad |I| \leq r$
are 
$L^{r}_{n}$-invariants with respect to the natural action (\ref{act}); they are
complete system of invariants in the sense of Weyl.
\end{itemize}

The whole formalism presented above can be realized in an arbitrary chart 
system 
$(V^{{\bf I},r},\psi^{{\bf I},r})$
on
${\rm Imm}T^{r}_{n}X$
so we have the central result:

\begin{thm}

The set
$
P^{r}_{n}X \equiv {\rm Imm} T_{n}^{r}X / L^{r}_{n}
$
has a unique differential manifold structure such that the canonical projection
$\rho^{r}_{n}$ 
is a submersion. The group action (\ref{action}) defines on
${\rm Imm} T_{n}^{r}X$
the structure of a right principal
$L_{n}^{r}$-bundle.
 
A chart system on
$
P^{r}_{n}X
$
adapted to this fibre bundle structure is formed from couples
$(W^{{\bf I},r},\Phi^{{\bf I},r})$
where:

\be
W^{{\bf I},r} = \left\{ j_{0}^{r}\gamma \in V^{r}| 
{\rm det}(x_{j}^{i_{k}}(j_{0}^{r}\gamma)) \not= 0 \right\}
\ee
and

\be
\Phi^{{\bf I},r} = (x^{i}_{I}, y^{\sigma}_{I}), \quad |I| \leq r.
\ee

In this case the local expression of the canonical projection is
$$
\rho^{r}_{n}(x^{i}_{I}, y^{\sigma}_{I}) = (x^{i},y^{\sigma}_{I}).
$$
\end{thm}

A point of
$P^{r}_{n}X$
containing a regular
$(r,n)$-velocity
$j^{r}_{0}\gamma$
is called an
$(r,n)$-{\it contact element}, or an
$r$-{\it contact element} of an
$n$-dimensional
submanifold of
$X$,
and is denoted by
$[j^{r}_{0}\gamma].$
As in the case of
$r$-jets, the point
$0 \in \R^{n}$
(resp.
$\gamma(0) \in X$)
is called the {\it source} (resp. the {\it target}) of 
$[j^{r}_{0}\gamma].$
The manifold
$P^{r}_{n}$
is called the
$(r,n)$-{\it Grassmannian bundle},
or simply a {\it higher order Grassmannian bundle} over
$X$.

Besides the quotient projection
$\rho^{r}_{n}: {\rm Imm} T^{r}_{n}X \to P^{r}_{n}$ 
we have for every 
$1 \leq s \leq r$,
the {\it canonical projection} of
$P^{r}_{n}X$
onto
$P^{s}_{n}X$
defined by
$\rho^{r,s}_{n}([j^{r}_{0}\gamma]) = [j^{s}_{0}\gamma]$ 
and the {\it canonical projection} of
$P^{r}_{n}X$
onto
$X$
defined by
$\rho^{r}_{n}([j^{r}_{0}\gamma]) = \gamma(0)$. 

On
$P_{n}^{r}X$
there are total differential operators; as expected, in the chart
$\rho^{r}_{n}(W^{{\bf I},r})$
they have the expression:

\be
\partial^{j_{1},\dots j_{k}}_{\sigma} \equiv {r_{1}! \dots r_{n}! \over k!}
{\partial \over \partial y^{\sigma}_{j_{1},\dots,j_{k}}}
\label{delta}
\ee

We note for further use the following formula:

\be
\partial^{i_{1},\dots,i_{k}}_{\sigma} \bar{y}^{\nu}_{j_{1},\dots,j_{k}}
= {\cal S}^{+}_{j_{1},\dots,j_{k}} P^{i_{1}}_{j_{1}} \dots P^{i_{k}}_{j_{k}}
Q^{\nu}_{\sigma}, \quad k = 1,...,r;
\label{Dy}
\ee
here we have defined:

\be
Q^{\sigma}_{\nu} \equiv \partial_{\nu} \bar{y}^{\sigma} - \bar{y}^{\sigma}_{i}
(\partial_{\nu} \bar{x}^{i}).
\label{P}
\ee

Next we define the {\it total derivative operators} on the Grassmann manifold:

\be
d_{i} \equiv {\partial \over \partial x^{i}} + 
\sum_{k=0}^{r-1} y^{\sigma}_{i,j_{1},\dots,j_{k}} 
\partial^{j_{1},\dots,j_{k}}_{\sigma} = 
{\partial \over \partial x^{i}} + 
\sum_{|J| \leq r-1} y^{\sigma}_{iJ} \partial^{J}_{\sigma}.
\label{D}
\ee

We note that:

\be
(\rho^{r}_{n})_{*} (z^{j}_{i} D_{j}) = d_{i}.
\label{Dd}
\ee

In particular, we have for any smooth function $f$ on 
$\rho^{r}_{n}(W^{r})$
the following formula:

\be
D_{i}(f \circ \rho^{r}_{n}) = x^{j}_{i} (d_{j}f) \circ \rho^{r}_{n}.
\label{dD}
\ee

The formula for the chart change on
$P^{r}_{n}X$.
can be written with this operators: let us consider two overlapping charts:
$(\rho^{r}_{n}(V^{r}),(x^{i},y^{\sigma}))$
and respectively 
$(\rho^{r}_{n}(\bar{V^{r}}),(\bar{x}^{i},\bar{y}^{\sigma}))$; 
then we have on the overlap:

\be
\bar{y}^{\sigma}_{iI} = P^{j}_{i} d_{j} \bar{y}^{\sigma}_{I}, 
\quad |I| \leq r-1 
\label{y-y-bar}
\ee
where
$P$
is the inverse of the matrix $Q$:

\be
Q^{l}_{p} \equiv d_{p}\bar{x}^{l}, \quad P^{j}_{i} Q^{l}_{j} = \delta^{l}_{i}.
\label{change}
\ee
We also note that:

\be
Q^{j}_{i} \bar{d}_{j} = d_{i}.
\label{DD}
\ee

\subsection{Contact Forms on Grassmann Manifolds}

By a {\it contact form} on
$P^{r}_{n}X$
we mean any form
$
\rho \in \Omega^{r}_{q}(PX)
$
verifying

\be
\left[j^{r}\gamma\right]^{*} \rho = 0
\ee
for any immersion
$\gamma: \R^{n} \rightarrow X$. 
We denote by
$\Omega^{r}_{q(c)}(PX)$
the set of contact forms of degree $q \leq n$. Here
$\left[ j^{r}\gamma\right] :\R^{n} \rightarrow P^{r}_{n}$
is given by:
$\left[ j^{r}\gamma\right] (t) \equiv \left[ j^{r}_{t}\gamma\right].$
We mention some of properties verified by these forms.

If one considers only the 
contact forms on an open set
$\rho^{r}_{n}(V^{r}) \subset P^{r}_{n}X$
then we emphasize this by writing
$\Omega^{r}_{q(c)}(V)$.
The ideal of all contact forms is denoted by
${\cal C}(\Omega^{r})$.
By elementary computations one finds out that, as in the case of a fibre
bundle, for any chart 
$(V,\psi)$
on $X$, every element of the set
$\Omega^{r}_{1(c)}(V)$
is a linear combination of the following expressions:

\be
\omega^{\sigma}_{j_{1},...,j_{k}} \equiv d y^{\sigma}_{j_{1},...,j_{k}} -
y^{\sigma}_{i,j_{1},...,j_{k}} d x^{i}, \quad k = 0,...,r-1
\label{o}
\ee
or, in multi-index notations:

\be
\omega^{\sigma}_{J} \equiv d y^{\sigma}_{J} - y^{\sigma}_{iJ} d x^{i}, \quad
|J| \leq r-1.
\label{o'}
\ee

We have the formula

\be
d \omega^{\sigma}_{J} = - \omega^{\sigma}_{Ji} \wedge  d x^{i}, \quad  
|J| \leq r-2.
\ee

Any form
$\rho \in \Omega^{r}_{q}(PX), \quad q = 2,...,n$
is contact {\it iff} it is generated by
$\omega^{\sigma}_{J}, \quad|J| \geq r-1$
and
$d\omega^{\sigma}_{I}, \quad |I| = r-1$.
In the end we present the transformation formula relevant for change of
charts. 

\begin{prop}
Let
$(V,\psi)$
and
$(\bar{V},\bar{\psi})$
two overlapping charts on $X$ and let
$(W^{r}, \Phi^{r}), \quad \Phi^{r} = (x^{i}, y^{\sigma}_{I}, x^{i}_{I})$
and
$(\bar{W}^{r}, \bar{\Phi}^{r}), \quad 
\bar{\Phi^{r}} = (\bar{x}^{i}, \bar{y}^{\sigma}_{I}, \bar{x}^{i}_{I})$
the corresponding charts on
$T^{r}_{n}X$.
Then the following formula is true on
$\rho^{r}_{n}(W^{r} \cap \bar{W}^{r}) \subset P^{r}_{n}X$:

\be
\bar{\omega}^{\sigma}_{I} = \sum_{|J|=1}^{|I|} (\partial^{J}_{\nu}
\bar{y}^{\sigma}_{I}) \omega^{\nu}_{J} - Q^{\sigma}_{I,\nu} \omega^{\nu}, 
\qquad 1 \leq |I| \leq r-1.
\label{transf}
\ee
where we have defined:
\be
Q^{\sigma}_{I,\nu} \equiv \partial_{\nu} \bar{y}^{\sigma}_{I} - 
\bar{y}^{\sigma}_{jI} (\partial_{\nu} \bar{x}^{j})
\label{Q_I}, \quad 0 \leq |I| \leq r-1
\ee
and 

\be
\bar{\omega}^{\sigma} = Q^{\sigma}_{\nu} \omega^{\nu}
\label{transf0}
\ee
where 
$Q^{\sigma}_{\nu}$
is given by the formula (\ref{P}).
\end{prop}

As a consequence we have:

\begin{cor}

If for a $q$-form has the expression

\be
\rho  = \sum_{p+s=k} \sum_{|J_{1}|,...,|J_{p}| \leq r-1}
\sum_{|I_{1}|=...=|I_{s}|=r-1} 
\omega^{\sigma_{1}}_{J_{1}} \cdots \wedge \omega^{\sigma_{p}}_{J_{p}} \wedge
d\omega^{\nu_{1}}_{I_{1}} \cdots \wedge d\omega^{\nu_{s}}_{I_{s}} 
\wedge \Phi^{J_{1},...,J_{p},I_{1},...,I_{s}}_{\sigma_{1},...,\sigma_{p},
\nu_{1},...,\nu_{s}}, \quad t \leq q
\label{contact-k}
\ee
is valid in one chart, then it is valid in any other chart.
\end{cor}

This corollary allows us to define for any
$q = 1,...,dim(J^{r}Y) = m{n+r \choose n}$  
a contact form with {\it order of contactness} $k$ to be any
$\rho \in \Omega^{r}_{q}$ 
such that it has in one chart (thereafter in any other chart) the expression
above. We denote these forms by
$\Omega^{r}_{q,k}$.
\newpage

\section{A Lagrange-Souriau Form on a Grassmann Manifold}

\subsection{Transformation Formulae and Invariant Conditions}

As in the preceding Section we consider a differential manifold $X$ and the
associated $(s,n)$-Grassmann manifold
$P^{s}_{n}X$.
We start we the following general result:

\begin{prop}
Let 
$\alpha \in \Omega^{s}_{q}(X)$
a $q$-differential form on
$P^{s}_{n}X$
verifying:

\be
i_{\xi} \alpha = 0
\ee
for any 
$\rho^{s,s-1}_{n}$-vertical vector field $\xi$ on
$P^{s}_{n}X$ 
(i.e. 
$(\rho^{s,s-1}_{n})_{*} \xi = 0$).
Then this form has the local expression:

\be
\alpha = \sum_{k=0}^{q} \sum_{|I_{1}|,\dots,|I_{k}| \leq s-1} 
T^{I_{1},\dots,I_{k}}_{\sigma_{1},\dots,\sigma_{k},i_{k+1},\dots,i_{q}}
\omega^{\sigma_{1}}_{I_{1}} \wedge \dots \wedge 
\omega^{\sigma_{k}}_{I_{k}} \wedge dx^{i_{k+1}} \wedge \dots \wedge dx^{i_{q}}
\label{alpha}
\ee
where
$T^{I_{1},\dots,I_{k}}_{\sigma_{1},\dots,\sigma_{k},i_{k+1},\dots,i_{q}}$
are smooth functions depending on the variables
$(x^{i},y^{\sigma},y^{\sigma}_{j},\dots,y^{\sigma}_{j_{1},\dots,j_{s}})$
and verify the (anti)symmetry property:

\be
T^{I_{P(1)},\dots,I_{P(k)}}_{\sigma_{P(1)},\dots,\sigma_{P(k)},
i_{Q(k+1)},\dots,i_{Q(q)}} = (-1)^{|P|+|Q|} 
T^{I_{1},\dots,I_{k}}_{\sigma_{1},\dots,\sigma_{k},i_{k+1},\dots,i_{q}};
\label{antisymmetry}
\ee
here $P$ is a permutation of the numbers 
$1,\dots,k$,
$Q$ is a permutation of the numbers
$k+1,\dots,q$
and
$|P|, |Q|$
are the signatures of these permutations. Moreover the following transformation
property is valid on the overlap of two charts:

\begin{eqnarray}
T^{I_{1},\dots,I_{p},\emptyset,\dots,\emptyset}_{\sigma_{1},\dots,\sigma_{p},
\dots,\sigma_{l},i_{l+1},\dots,i_{q}} = 
\sum_{k=p}^{l} {k \choose p} {q-k \choose q-l} 
\sum_{|J_{1}| \geq |I_{1}|} \cdots \sum_{|J_{p}| \geq |I_{p}|} 
\sum_{|J_{p+1}|,\dots, |J_{k}| \leq s-1} 
(\partial^{I_{1}}_{\sigma_{1}} \bar{y}^{\nu_{1}}_{J_{1}}) \dots
(\partial^{I_{p}}_{\sigma_{p}} \bar{y}^{\nu_{p}}_{J_{p}})
\nonumber \\ 
Q^{\nu_{p+1}}_{J_{p+1},\sigma_{p+1}} \cdots Q^{\nu_{k}}_{J_{k},\sigma_{k}} 
(\partial_{\sigma_{k+1}} \bar{x}^{j_{k+1}}) \cdots
(\partial_{\sigma_{l}} \bar{x}^{j_{l}}) 
Q^{j_{l+1}}_{i_{l+1}} \cdots Q^{j_{q}}_{i_{q}}
\bar{T}^{J_{1},\dots,J_{k}}_{\nu_{1},\dots,\nu_{k},j_{k+1},\dots,j_{q}}
\label{TT-bar}
\end{eqnarray}
where 
$p \leq l$,
$I_{1},\dots,I_{p} \not= \emptyset$
and the notations
$Q^{\nu}_{J,\sigma}$
and
$Q^{j}_{i}$
have been introduced according to the formulae (\ref{Q_I}) and (\ref{change})
respectively. 
\end{prop}

{\bf Proof:} Follows by elementary computations from (\ref{transf}) and 
(\ref{transf0}).
\qed

We denote the space of these forms by
$\Omega^{s}_{q,\xi}(X)$.
As a a consequence of this proposition we have some corollaries.

\begin{cor}

Let
$\alpha \in \Omega^{s}_{q,\xi}(X)$
having the local expression (\ref{alpha}). Then the following relations have a
intrinsic global meaning:

\be
T^{I_{1},\dots,I_{k}}_{\sigma_{1},\dots,\sigma_{k},i_{k+1},\dots,i_{q}} = 0,
\quad |I_{1}| + \cdots + |I_{k}| \geq t
\label{T1}
\ee
for any $t \in \N$.

\end{cor}

We denote the subset of these forms by
$\Omega^{s,t}_{q,\xi}(X)$.

\begin{cor}

Let
$\alpha \in \Omega^{s}_{q,\xi}(X)$
having the local expression (\ref{alpha}). Then the following relations have a
intrinsic global meaning:

\be
T^{I_{1}}_{\sigma_{1},i_{2},\dots,i_{q}} = 0, \quad 
\forall I_{1} \not= \emptyset.
\label{T2}
\ee
\end{cor}

If the condition (\ref{T2}) is fulfilled we say that the form $\alpha$ verifies
the {\it Lepage condition}.  We denote the subset of the forms verifying the
conditions (\ref{T1}) and (\ref{T2}) by 
$\Omega^{s,t,Lep}_{q,\xi}(X)$.

Similarly we have

\begin{cor}

Let
$\alpha \in \Omega^{s}_{q,\xi}(X)$
having the local expression (\ref{alpha}). Then the following relations have a
intrinsic global meaning:

\be
T_{i_{1},\dots,i_{q}} = 0.
\label{T0}
\ee
\end{cor}

We make another useful notations:
$$
\Omega^{s,t}_{q,\xi,k}(X) \equiv \Omega^{s}_{q,k}(X) \cap
\Omega^{s,t}_{q,\xi}(X);
$$
these are contact forms with the order of contactness equal to $k$ and the
definition is globally true.

The proof of these corollaries are elementary following from the proposition
above. 

A little more complicated is a result for which we need the tensorial notations
introduced in \cite{Gr3}.  We introduce the tensor spaces: 
$$ 
{\cal H}_{k} \equiv {\cal F}^{(-)}(\R^{n}) \otimes
\underbrace{{\cal F}^{(+)}(\R^{n}) \otimes \cdots \otimes
{\cal F}^{(+)}(\R^{n})}_{k-times}
$$
where
${\cal F}^{(\pm)}(\R^{n})$
are the symmetric (coresp. $+$) and the antisymmetric (coresp. $-$) Fock
spaces. We have the well known decomposition in subspaces with fixed number of
``bosons" and ``fermions": 
$$ 
{\cal H}_{k} = \oplus_{l=0}^{s} \oplus_{r_{1},...,r_{k} \geq 0} 
{\cal H}_{l,r_{1},...,r_{k}}.  
$$

We make the convention that
${\cal H}_{k,l_{1},...,l_{k}} \equiv 0$
if any one of the indices
$l,r_{1},...,r_{k}$
is negative or if
$l > n$.
Then we can consider
$T^{I_{1},...,I_{k}}_{\sigma_{1},...,\sigma_{k},i_{k+1},...,i_{q}}$
as the components of a tensor
$
T_{\sigma_{1},...,\sigma_{k}} \in {\cal H}_{q-k,|I_{1}|,...,|I_{k}|}.
$  

We can write in an compact way the next corollary if we use the creation and
the annihilation fermionic operators 
$a^{*i}, a_{i}, \quad (i = 1,...,n)$ 
and
the corresponding creation and annihilation bosonic operators
$b^{*}_{(p)i}, b^{i}_{(p)} \quad (p = 1,...,k; i = 1,...,n).$ 
We are using conventions somewhat different from that used in quantum
mechanics (see \cite{Gr3}) which amount to a rescalation of the usual
expressions in every subspace with fixed ``number of particles". Finally, one 
introduces the operators
$B_{p}, \quad (p = 1,...,k)$
according to:

\be
B_{p}\vert_{{\cal H}_{k}} \equiv (-1)^{k} b^{*}_{(p)i} a^{*i}, .
\label{B}
\ee

Now we can formulate

\begin{cor}
Let 
$\alpha \in \Omega^{s,t}_{q,\xi}(X)$.
Then the conditions:

\be
(B^{*}_{p}T)^{I_{1},...,I_{k}}_{\sigma_{1},...,\sigma_{k},i_{k+2},...,i_{q}} 
= 0, \quad |I_{1}| + \cdots + |I_{k}| = t-2, 
\quad p = 1,\dots,k, \quad k =1,\dots,q
\label{T3}
\ee
are globally defined.
\end{cor}

{\bf Proof:} Because of the antisymmetry property (\ref{antisymmetry}) it is
sufficient to consider only the case 
$p = 1$.
If we use the definition of the adjoint
$B^{*}_{1}$
and the transformation rule given in the preceding proposition we obtain:
\begin{eqnarray}
(B^{*}_{1}T)^{I_{1},...,I_{k}}_{\sigma_{1},...,\sigma_{k},i_{k+2},...,i_{q}}
= {\rm const} \times
T^{lI_{1},...,I_{k}}_{\sigma_{1},...,\sigma_{k},l,i_{k+2},...,i_{q}} =
\nonumber \\
{\rm const} \times \sum_{|J_{1}| = |I_{1}|+1}^{s} \sum_{|J_{2}| = |I_{2}|}^{s} 
\cdots \sum_{|J_{k}| = |I_{k}|}^{s} (\partial^{lI_{1}}_{\sigma_{1}} 
\bar{y}^{\nu_{1}}_{J_{1}}) Q^{i}_{l}
\cdots \bar{T}^{J_{1},...,J_{k}}_{\nu_{1},...,\nu_{k},i,j_{k+2},...,j_{q}}.
\nonumber
\end{eqnarray}

By $\dots$ we mean here the other factors from the formula (\ref{TT-bar}) for
which the full expression is not needed.
But it is not very hard to see that the conditions (\ref{T1}) and (\ref{T3}) 
which are verified by the form $\alpha$ impose:
$|J_{1}| = |I_{1}|+1, |J_{2}| = |I_{2}|, \dots, |J_{k}| = |I_{k}|$.
Let us take therefore
$|I_{1}| = p$
and
$J_{1} \equiv \{u_{0},\dots,u_{p}\}$;
we obtain from (\ref{Dy}) that
$$
Q^{i}_{l} (\partial^{lI_{1}}_{\sigma_{1}} \bar{y}^{\nu_{1}}_{J_{1}}) =
S^{+}_{u_{0},\dots,u_{p}} \delta^{i}_{u_{0}} \partial^{I_{1}}_{\sigma_{1}}
\bar{y}^{\nu_{1}}_{u_{1},\dots,u_{p}}.
$$

So we get finally
\begin{eqnarray}
(B^{*}_{1}T)^{I_{1},...,I_{k}}_{\sigma_{1},...,\sigma_{k},i_{k+2},...,i_{q}} =
{\rm const} \times 
\sum_{|J_{1}| = |I_{1}|}^{s} \cdots \sum_{|J_{k}| = |I_{k}|}^{s}
(\partial^{I_{1}}_{\sigma_{1}} \bar{y}^{\nu_{1}}_{J_{1}}) 
\cdots \bar{T}^{iJ_{1},...,J_{k}}_{\nu_{1},...,\nu_{k},i,j_{k+2},...,j_{q}} =
\nonumber \\
\sum_{|J_{1}| = |I_{1}|}^{s} \cdots \sum_{|J_{k}| = |I_{k}|}^{s}
(\partial^{I_{1}}_{\sigma_{1}} \bar{y}^{\nu_{1}}_{J_{1}}) 
\cdots (B^{*}_{1}\bar{T})^{J_{1},...,J_{k}}_{\nu_{1},...,\nu_{k},
j_{k+2},...,j_{q}}
\end{eqnarray}
and the assertion from the statement follows.
\qed

\begin{rem}
The conditions (\ref{T3}) mean that the tensors entering in the expression of
the form $\alpha$ are traceless. 
\end{rem}

We denote the subset of the forms verifying the conditions (\ref{T1}),
(\ref{T2}) (i.e Lepage) and (\ref{T3}) (i.e. tracelessness) by
$\Omega^{s,t,Lep}_{q,\xi,tr}(X)$. 
We finally stress again that all the spaces of the type
$\Omega^{\dots}_{\dots}$
are globally defined.
\newpage

\subsection{The Definition of the Lagrange-Souriau Form}

We start this Subsection with a general result. First we introduce some
notations which will be very useful in the following. We define
$$ 
{\cal H}_{l,k} \equiv \oplus_{r_{1},...,r_{k}=0}^{s-1}  
{\cal H}_{l,r_{1},...,r_{k}}.  
$$

Let then
$
T \equiv \{T_{k}\}; \quad T_{k} \in {\cal H}_{q-k,k}
$
be an ensemble of tensors as the ones appearing in the structure of the form
(\ref{alpha}). We define some new tensors as follows:

- $\delta_{1}: {\cal H}_{q-k,k} \rightarrow {\cal H}_{q-k,k+1}$
by:

\begin{eqnarray}
(\delta_{1} T)^{I_{1},\dots,I_{k}}_{\sigma_{1},\dots,\sigma_{k},
i_{k+1},\dots,i_{q+1}} \equiv 
S^{-}_{(I_{1},\sigma_{1}),\dots,(I_{k},\sigma_{k})}
\partial^{I_{1}}_{\sigma_{1}} 
T^{I_{2},\dots,I_{k}}_{\sigma_{2},\dots,\sigma_{k},i_{k+1},\dots,i_{q+1}} =
\nonumber \\
{1 \over k} \sum_{p=1}^{k} (-1)^{p-1} \partial^{I_{p}}_{\sigma_{p}} 
T^{I_{1},\dots,\hat{I}_{p},\dots,I_{k}}_{\sigma_{1},\dots,\hat{\sigma}_{p},
\dots,\sigma_{k},i_{k+1},\dots,i_{q+1}};
\label{delta1}
\end{eqnarray}

- $\delta_{2}: {\cal H}_{q-k,k} \rightarrow {\cal H}_{q-k+1,k}$
by:

\begin{eqnarray}
(\delta_{2} T)^{I_{1},\dots,I_{k}}_{\sigma_{1},\dots,\sigma_{k},
i_{k+1},\dots,i_{q+1}} \equiv (-1)^{q}
S^{-}_{i_{k+1},\dots,i_{q+1}} d_{i_{q+1}}
T^{I_{1},\dots,I_{k}}_{\sigma_{1},\dots,\sigma_{k},i_{k+1},\dots,i_{q}} +
\sum_{p=1}^{k} (B_{p}T)^{I_{1},\dots,I_{k}}_{\sigma_{1},\dots,\sigma_{k},
i_{k+1},\dots,i_{q+1}} \nonumber \\
= {1 \over q-k} \sum_{p=k+1}^{q+1} (-1)^{p-1} d_{i_{p}}
T^{I_{1},\dots,I_{k}}_{\sigma_{1},\dots,\sigma_{k},i_{k+1},\dots,
\hat{i_{p}},\dots,i_{q+1}} +
\sum_{p=1}^{k} (B_{p}T)^{I_{1},\dots,I_{k}}_{\sigma_{1},\dots,\sigma_{k},
i_{k+1},\dots,i_{q+1}};
\label{delta2}
\end{eqnarray}
and

\be
\delta T \equiv \delta_{1} T + \delta_{2} T
\label{delta12}
\ee

Then we have:

\begin{prop}

Let 
$\alpha \in \Omega^{s,s-1}_{q,\xi}(X)$
be arbitrary. If we put it in the local coordinates according to (\ref{alpha})
then it follows that its exterior differential is given by:

\begin{eqnarray}
d\alpha = \sum_{k=0}^{q} \sum_{|I_{0}|=s} 
\sum_{|I_{1}|,\dots,|I_{k}| \leq s-1} \left[ \partial^{I_{0}}_{\sigma_{0}} 
T^{I_{1},\dots,I_{k}}_{\sigma_{1},\dots,\sigma_{k},i_{k+1},\dots,i_{q}} +(k+1) 
(B_{0}T)^{I_{0},\dots,I_{k}}_{\sigma_{0},\dots,\sigma_{k},i_{k+1},\dots,i_{q}}
\right] \nonumber \\
dy^{\sigma_{0}}_{I_{0}} \wedge \omega^{\sigma_{1}}_{I_{1}} \wedge 
\dots \wedge \omega^{\sigma_{k}}_{I_{k}} \wedge dx^{i_{k+1}} \wedge \dots 
\wedge dx^{i_{q}} \nonumber \\
+ \sum_{k=0}^{q+1} \sum_{|I_{1}|,\dots,|I_{k}| \leq s-1} 
(\delta T)^{I_{1},\dots,I_{k}}_{\sigma_{1},\dots,\sigma_{k},
i_{k+1},\dots,i_{q+1}}
\omega^{\sigma_{1}}_{I_{1}} \wedge \dots \wedge \omega^{\sigma_{k}}_{I_{k}} 
\wedge dx^{i_{k+1}} \wedge \dots \wedge dx^{i_{q+1}}.
\label{dalpha}
\end{eqnarray}
\label{alpha-dalpha}
\end{prop}

The proof follows by elementary, but somewhat tedious computations. We also
have 

\begin{cor}

The following formula is true:

\be
\delta^{2} = 0.
\ee

As a consequence we also have:

\be
\delta_{1}^{2} = 0, \quad 
\delta_{1} \delta_{2} + \delta_{2} \delta_{1} = 0, \quad
\delta_{2}^{2} = 0.
\ee
\label{delta-square}
\end{cor}

{\bf Proof:} We write the formula (\ref{dalpha}) as follows:
$$
d\alpha = \sum_{k=0}^{q+1} \sum_{|I_{1}|,\dots,|I_{k}| \leq s-1} 
(\delta T)^{I_{1},\dots,I_{k}}_{\sigma_{1},\dots,\sigma_{k},
i_{k+1},\dots,i_{q+1}}
\omega^{\sigma_{1}}_{I_{1}} \wedge \dots \wedge \omega^{\sigma_{k}}_{I_{k}} 
\wedge dx^{i_{k+1}} \wedge \dots \wedge dx^{i_{q+1}} + \cdots
$$
where by $\cdots$ we mean terms containing at least one differential
$dy^{\sigma}_{I}, \quad |I| = s$.

If we use now $d^{2} = 0$ and iterate the formula above we get
$$
\sum_{k=0}^{q+2} \sum_{|I_{1}|,\dots,|I_{k}| \leq s-1} 
(\delta^{2} T)^{I_{1},\dots,I_{k}}_{\sigma_{1},\dots,\sigma_{k},
i_{k+1},\dots,i_{q+2}}
\omega^{\sigma_{1}}_{I_{1}} \wedge \dots \wedge \omega^{\sigma_{k}}_{I_{k}} 
\wedge dx^{i_{k+1}} \wedge \dots \wedge dx^{i_{q+2}} + \cdots = 0
$$
which implies
$$
\delta^{2} T_{k} = 0 \Longleftrightarrow
\left[ \delta_{1}^{2} + (\delta_{1} \delta_{2} + \delta_{2} \delta_{1}) + 
\delta_{2}^{2} \right] T_{k} = 0, \quad k = 0,\dots,q 
$$
 
But the three terms from the left hand side belong to the subspaces:
$
{\cal H}_{q-k,k+2}, {\cal H}_{q-k+1,k+1}$
and
$
{\cal H}_{q-k+2,k}
$
respectively, so they must be zero separately. 
$\qed$

\begin{rem}
These formulae can be obtained directly from the definitions (\ref{delta1})
and (\ref{delta2}). We also note that the operators 
$\delta, \delta_{1}, \delta_{2}$
verify a BRST type algebra and can be used to build a bicomplex in the space of
the tensors 
$\{T_{k}\}$
(see \cite{An}). However, these definitions are only local, i.e. they do not
have an intrinsic geometrical meaning.
\end{rem}

As a corollary of the preceding proposition we have:

\begin{cor}
If
$\alpha \in \Omega^{s,s-1}_{q,\xi}(X)$
is closed, then the following relations are true:

\begin{eqnarray}
\partial^{I_{0}}_{\sigma_{0}} 
T^{I_{1},\dots,I_{k}}_{\sigma_{1},\dots,\sigma_{k},i_{k+1},\dots,i_{q}} +(k+1) 
(B_{0}T)^{I_{0},\dots,I_{k}}_{\sigma_{0},\dots,\sigma_{k},i_{k+1},\dots,i_{q}}
= 0, \nonumber \\ \quad |I_{0}| = s, \quad |I_{1}|,\dots,|I_{k}| \leq s-1, 
\quad k = 0,\dots,q
\label{Ts}
\end{eqnarray}
and

\be
(\delta T)^{I_{1},\dots,I_{k}}_{\sigma_{1},\dots,\sigma_{k},
i_{k+1},\dots,i_{q+1}} = 0, \quad |I_{1}|,\dots,|I_{k}| \leq s-1, 
\quad k = 0,\dots,q+1.
\label{deltaT}
\ee

In particular we have

\be
\partial^{I_{0}}_{\sigma_{0}} 
T^{I_{1},\dots,I_{k}}_{\sigma_{1},\dots,\sigma_{k},i_{k+1},\dots,i_{q}} = 0, 
\quad |I_{0}| = s, \quad |I_{1}| + \dots + |I_{k}| \geq 1, 
\quad k = 0,\dots,q
\ee
i.e. the tensors
$
T^{I_{1},\dots,I_{k}}_{\sigma_{1},\dots,\sigma_{k},i_{k+1},\dots,i_{q}} 
\quad |I_{0}| = s, \quad |I_{1}| + \dots + |I_{k}| \geq 1, \quad k = 0,\dots,q
$
depend only on the variables
$(x^{i},y^{\sigma},y^{\sigma}_{j},\dots,y^{\sigma}_{j_{1},\dots,j_{s-1}})$
and 
$
T^{\emptyset,\dots,\emptyset}_{\sigma_{1},\dots,\sigma_{k},
i_{k+1},\dots,i_{q}}, \quad k = 0,\dots,q
$
depend on 
$y^{\sigma}_{j_{1},\dots,j_{s}}$
only linearly.
\label{lin}
\end{cor}

\begin{rem}
We anticipate somewhat and draw the attention to the extremely interesting
equation (\ref{Ts}); we will prove that if the tensors
$\{T_{k}\}$
are verifying this equation, they depend on the highest-order derivatives
through some polynomial expressions called hyper-Jacobians \cite{BCO},
\cite{Ol1}. 
\end{rem}

This simple structure of the tensors leads to

\begin{prop}

Let
$\alpha \in \Omega^{s,s-1}_{q,\xi}(X)$
be closed. Then $\alpha$ is
$\rho^{s,s-1}_{n}$-projectable,
i.e. there exists a $q$-form 
$\alpha_{0} \in \Omega^{s-1}_{q}(X)$
such that

\be
\alpha = (\rho^{s,s-1}_{n})^{*} \alpha_{0}.
\label{alpha0}
\ee

Moreover, the form $\alpha_{0}$ is closed.
\label{project}
\end{prop}

{\bf Proof:}
We exhibit the dependence of the form $\alpha$ on the highest-order
derivatives; according to the preceding corollary these derivatives can appear
in two places: in the coefficients 
$
T^{\emptyset,\dots,\emptyset}_{\sigma_{1},\dots,\sigma_{k},
i_{k+1},\dots,i_{q}}, \quad k = 0,\dots,q
$
and in the contact forms
$\omega^{\sigma}_{I}, \quad |I| = s-1.$
It is not very hard to write now $\alpha$ as follows:
\begin{eqnarray}
\alpha = \sum_{k=0}^{q} 
T^{\emptyset,\dots,\emptyset}_{\sigma_{1},\dots,\sigma_{k},i_{k+1},\dots,i_{q}}
\omega^{\sigma_{1}} \wedge \dots \wedge \omega^{\sigma_{k}} 
\wedge dx^{i_{k+1}} \wedge \dots \wedge dx^{i_{q}} + \nonumber \\
\sum_{k=1}^{q} \sum_{|I_{1}| = s-1} k 
T^{I_{1},\emptyset,\dots,\emptyset}_{\sigma_{1},\dots,\sigma_{k},
i_{k+1},\dots,i_{q}}
\omega^{\sigma_{1}}_{I_{1}} \wedge \omega^{\sigma_{2}} \dots \wedge 
\omega^{\sigma_{k}} \wedge dx^{i_{k+1}} \wedge \dots \wedge dx^{i_{q}} +
\alpha'\nonumber
\end{eqnarray}
where the form $\alpha'$ is 
$\rho^{s,s-1}_{n}$-projectable.

The expression above is, as said before, at most linear in the highest-order
derivatives. One computes explicitly the coefficient of this derivatives and
finds out that they are zero according to (\ref{Ts}). This proves the first
assertion. The closedness of the form $\alpha_{0}$ follows from the
surjectivity of the map
$\rho^{s,s-1}_{n}.$
$\qed$

It is clear that there are strong conditions on any closed form 
$\alpha \in \Omega^{s,s-1,Lep}_{q,\xi,tr}(X)$.
We will give a structure theorem for such a form in the case
$q = n+1$ 
which is relevant for physical applications. First we note that in this case we
have

\be
T^{\emptyset}_{\sigma,i_{1},\dots,i_{n}} = T_{\sigma} 
\varepsilon_{i_{1},\dots,i_{n}} 
\label{T-sigma}
\ee
for some smooth functions 
$T_{\sigma}$
on
$P^{s}_{n}X.$
Here
$\varepsilon_{i_{1},\dots,i_{n}}$
is the completely antisymmetric tensor.

\begin{thm}

Let
$\alpha \in \Omega^{s,s,Lep}_{n+1,\xi,tr}(X)$
be closed. Then $\alpha$ admits the following decomposition:

\be
\alpha = T_{0} + dT_{1}
\label{decomp}
\ee
where:
\begin{itemize}
\item
$T_{1} \in \Omega^{s,s-2}_{n,\xi,2}(X)$
and
$dT_{1} \in \Omega^{s,s-1}_{n+1,\xi,2}(X)$.
\item
$T_{0} \in \Omega^{s,s-1,1}_{n+1,\xi}(X)$
has the local structure given by the formula (\ref{alpha}) with the tensors
$T_{k}$
given by formulae of the type

\be
T_{k} = P_{k} T_{\sigma}, \quad k = 2,\dots,n+1
\label{T-polyn}
\ee
with
$P_{k}$
some linear differential operators which can be reccurrsively determined. 
\end{itemize}
\label{canonical}
\end{thm}

{\bf Proof:}
The proof relies heavily on induction. 

(i) We first consider the equation (\ref{Ts}) for 
$k = 1$.
If
$I_{1} \not= \emptyset$
then, because of (\ref{T1}) and (\ref{T2}) an identity is obtained. If
$I_{1} = \emptyset$
then we get

\be
\partial^{I_{0}}_{\sigma_{0}} 
T^{\emptyset}_{\sigma_{1},i_{2},\dots,i_{n+1}} + 2 
(B_{0}T)^{I_{0},\emptyset}_{\sigma_{0},\sigma_{1},i_{2},\dots,i_{n+1}} = 0,
\quad |I_{0}| = s
\ee
or in compact notations:

\be
B_{0} T_{2} = - C_{2}
\label{TC}
\ee
where the components of the tensor
$C_{2}$
are:

\be
C^{I_{0},\emptyset}_{\sigma_{0},\sigma_{1},i_{2},\dots,i_{n+1}} \equiv
{1 \over 2} \partial^{I_{0}}_{\sigma_{0}} 
T^{\emptyset}_{\sigma_{1},i_{2},\dots,i_{n+1}}.
\ee

It is not very hard to prove that the previous relation is equivalent to the
result obtained after applying the operator
$B^{*}_{0}$
to it. Indeed, one knows (see \cite{Gr3}) that 
$
B^{*}_{0}B_{0}
$
and
$
B_{0}B^{*}_{0}
$
are (up to some constants) projector operators. So, the equation above is
equivalent to two relations: one obtained by applying the first operator and the
second by applying the second operator. But one sees easily that in the first
case an identity is obtained (because
$B_{0} {\cal H}_{n,l} \subset {\cal H}_{n+1,l+1} = \{0\}$)
so one obtains that the relation above, which is of the type
$X = 0$
is equivalent to
$B_{0} B^{*}_{0} X = 0$.
Taking the scalar product with $X$ one obtains from here
$\Vert B^{*}_{0} X \Vert^{2} = 0 \Longrightarrow B^{*}_{0} X = 0.$ 

So we have equivalently with (\ref{TC}):

$$
B^{*}_{0} B_{0} T_{2} = - B^{*}_{0} C_{2}.
$$
 
Let 
$N_{0}$
and
$N_{f}$
be the bosonic (resp. fermionic) number operators; one knows that (\cite{Gr3})
$$
\{B^{*}_{0}B_{0} \} = N_{0} - N_{f} + n {\bf 1};
$$

We use this relation in the equation above and take into account the
tracelessness condition (\ref{T3}). One obtains  
$$ 
s T_{2} = -  B^{*}_{0} C_{2}
$$ 
or, with full index notations:

\be
T^{I_{0},\emptyset}_{\sigma_{1},\sigma_{2},i_{2},\dots,i_{n+1}} = 
- {1 \over 2} n \partial^{I_{1}i_{1}}_{\sigma_{1}} 
T^{\emptyset}_{\sigma_{2},i_{1},\dots,i_{n+1}}, \quad |I_{1}| = s-1.
\label{T^s-1}
\ee

(ii) It is convenient to introduce the lexicographic order relation on the set
of $k$-uples
$
(I_{1},\dots,I_{k}) 
$
according to:
$$
(I_{1},\dots,I_{k}) > (J_{1},\dots,J_{k}) \quad {\it iff} \quad
|I_{1}| > |J_{1}|, \quad {\rm or} \quad
|I_{1}| = |J_{1}| \quad {\rm and} \quad |I_{2}| > |J_{2}|, \quad {\rm etc.}
$$

We can prove now by induction that we have in general for 
$k = 0, \dots, n+1:$

\be
T^{I_{1},\dots,I_{k}}_{\sigma_{1},\dots,\sigma_{k},i_{k+1},\dots,i_{n+1}} =
P^{I_{1},\dots,I_{k}}_{\sigma_{1},\dots,\sigma_{k},i_{k+1},\dots,i_{n+1}}
(T_{\nu}) + 
(\delta \Lambda)^{I_{1},\dots,I_{k}}_{\sigma_{1},\dots,\sigma_{k},
i_{k+1},\dots,i_{n+1}}.
\label{T_k}
\ee
Here
$
P^{I_{1},\dots,I_{k}}_{\sigma_{1},\dots,\sigma_{k},i_{k+1},\dots,i_{n+1}}
$
are linear differential operators verifying conditions of the type (\ref{T1}),
(\ref{T2}) and (\ref{T3}) and can be determined reccurrsively; also
$
\Lambda^{I_{1},\dots,I_{k}}_{\sigma_{1},\dots,\sigma_{k},i_{k+1},\dots,i_{n}}
$
are tensors depending on the variables
$
(x^{i},y^{\sigma},y^{\sigma}_{j},\dots,y^{\sigma}_{j_{1},\dots,j_{s-1}})
$
and verifying the following formulae

\be
\Lambda^{I_{1},\dots,I_{k}}_{\sigma_{1},\dots,\sigma_{k},i_{k+1},\dots,i_{n}} =
0, \quad |I_{1}| + \dots + |I_{k}| \geq s - 2
\label{lambda1}
\ee
and

\be
(\delta\Lambda)^{I_{1},\dots,I_{k}}_{\sigma_{1},\dots,\sigma_{k},
i_{k+1},\dots,i_{n+1}} = 0, \quad |I_{1} + \dots + |I_{k}| \geq s - 1.
\label{lambda2}
\ee

We will use a double induction procedure. More explicitly we prove the relation
(\ref{T_k}) by induction on (increasing values of) $k$ and for every fixed $k$
we prove the relation by induction on the (decreasing) lexicographic ordering.

(iii) We consider more explicitly the case $k = 2$. In this case we have been
able to find the explicit expression for the polynomials 
$P^{\dots}_{\dots}$. 
Namely we prove that (\ref{T_k}) has the following explicit form for the case
$k = 2$:

\be
T^{I_{1},I_{2}}_{\sigma_{1},\sigma_{2},i_{2},\dots,i_{n}} =
- {1 \over 2} n \sum_{|J| \leq s-1-|I_{1}|-|I_{2}|} (-1)^{|J|+|I_{2}|}
{|J|+|I_{2}| \choose |J|} d_{J} \partial^{I_{1}I_{2}Ji_{1}}_{\sigma_{1}}
T^{\emptyset}_{\sigma_{2},i_{1},\dots,i_{n}} + 
(\delta \Lambda)^{I_{1},I_{2}}_{\sigma_{1},\sigma_{2},i_{2},\dots,i_{n}}
\label{T_2}
\ee
where we have (\ref{lambda1}) and (\ref{lambda2}):

\be
\Lambda^{I_{1},I_{2}}_{\sigma_{1},\sigma_{2},i_{3},\dots,i_{n}} =
0, \quad |I_{1}| + |I_{2}| \geq s - 2, 
\quad
\Lambda^{I_{1}}_{\sigma_{1},i_{2},\dots,i_{n}} = 0, \quad \forall |I_{1}| \leq 
s - 1.
\ee

The idea is to use the equation (\ref{deltaT}) for $k = 2$:

\begin{eqnarray}
{1\over 2} \left(
\partial^{I_{1}}_{\sigma_{1}} T^{I_{2}}_{\sigma_{2},i_{3},\dots,i_{n+2}} -
\partial^{I_{2}}_{\sigma_{2}} T^{I_{1}}_{\sigma_{1},i_{3},\dots,i_{n+2}}
\right) +
{1 \over n} \sum_{p=3}^{n+2} (-1)^{p-1}  d_{i_{p}}
T^{I_{1},I_{2}}_{\sigma_{1},\sigma_{2},i_{3},\dots,\hat{i_{p}},\dots,i_{n+2}} +
\nonumber \\
(B_{1}T)^{I_{1},I_{2}}_{\sigma_{1},\sigma_{2},i_{3},\dots,i_{n+2}} +
(B_{2}T)^{I_{1},I_{2}}_{\sigma_{1},\sigma_{2},i_{3},\dots,i_{n+2}} = 0,
\quad |I_{1}|, |I_{2}| \leq s - 1.
\label{deltaT2}
\end{eqnarray}

One can write this equation in the equivalent form:

\be
(B_{1}T)^{I_{1},I_{2}}_{\sigma_{0},\sigma_{1},i_{2},\dots,i_{n+1}} = 
- C^{I_{1},I_{2}}_{\sigma_{0},\sigma_{1},i_{2},\dots,i_{n+1}};
\label{TC2}
\ee
because 
$I_{1} \not= \emptyset$
we only need:

\be
C^{I_{1},I_{2}}_{\sigma_{0},\sigma_{1},i_{2},\dots,i_{n+1}} \equiv
{1 \over 2} 
\partial^{I_{1}}_{\sigma_{1}} T^{I_{2}}_{\sigma_{2},i_{3},\dots,i_{n+1}} +
{1 \over n} \sum_{p=2}^{n+1} (-1)^{p}  d_{i_{p}}
T^{I_{1},I_{2}}_{\sigma_{1},\sigma_{2},i_{2},\dots,\hat{i_{p}},\dots,i_{n+1}}
+ (B_{2}T)^{I_{1},I_{2}}_{\sigma_{1},\sigma_{2},i_{2},\dots,i_{n+1}}.
\ee

With the same argument as at (i) we see that the equation (\ref{TC2}) above is
equivalent with the result obtained by applying to it the operator 
$B^{*}_{1}.$
We note that the tensors appearing in the right hand side of this relation have
the lexicographic order strictly greater the the lexicographic order of the
similar tensors from the left hand side. 
Now one can start by induction from the maximal lexicographic order (which
corresponds to  
$T^{I_{1},\emptyset}_{\dots}$
for $|I_{1}| = s - 1$
and it is known according to (\ref{T^s-1})) and prove formula (\ref{T_2}).
The computations are extremely long and tedious but straightforward and are
omitted. 

(iv) In the final step one proves (\ref{T_k}) assuming that it has been already
proved for 
$2, 3,\dots,k - 1.$
For this, one proceeds as above. First one writes (\ref{deltaT}) in the form
(\ref{TC2}) i.e.

\be 
(B_{1}T)^{I_{1},\dots,I_{k}}_{\sigma_{1},\dots,\sigma_{k},
i_{k+1},\dots,i_{n+2}} = - C^{I_{1},\dots,I_{k}}_{\sigma_{1},\dots,\sigma_{k},
i_{k+1},\dots,i_{n+2}}
\ee
where:

\begin{eqnarray}
C^{I_{1},\dots,I_{k}}_{\sigma_{1},\dots,\sigma_{k},
i_{k+1},\dots,i_{n+2}} \equiv 
- {1 \over k} \sum_{p=1}^{k} (-1)^{p-1} \partial^{I_{p}}_{\sigma_{p}} 
T^{I_{1},\dots,\hat{I}_{p},\dots,I_{k}}_{\sigma_{1},\dots,\hat{\sigma}_{p},
\dots,\sigma_{k},i_{k+1},\dots,i_{n+2}} +
\nonumber \\
{1 \over n-k+2} \sum_{p=k+1}^{n+2} (-1)^{p-1} d_{i_{p}}
T^{I_{1},\dots,I_{k}}_{\sigma_{1},\dots,\sigma_{k},i_{k+1},\dots,
\hat{i_{p}},\dots,i_{n+2}} +
\sum_{p=1}^{k} (B_{p}T)^{I_{1},\dots,I_{k}}_{\sigma_{1},\dots,\sigma_{k},
i_{k+1},\dots,i_{n+2}}.
\label{TCk}
\end{eqnarray}

One makes the same observations as in (iii), namely that the equation
(\ref{TCk}) above is equivalent with the result obtained by applying to it the
operator 
$B^{*}_{1}$ 
and that the tensors appearing in the right hand side of this relation have the
lexicographic order strictly greater the the lexicographic order of the similar
tensors from the left hand side or are already known from the induction
hypothesis. So it remains to prove the formula (\ref{T_k}) for the highest
lexicographical order which can be done as in (i). The induction argument is,
as in (iii), based on rather complicated computations which are skipped.

The induction is finished and the formula (\ref{T_k}) is proved. It remains to
match it with the formula (\ref{decomp}). Indeed, the first expression leads to
the first component $T_{0}$ of $\alpha$ and it remains to apply (\ref{dalpha})
to prove that the second expression leads to the second component $dT_{1}.$
Moreover, the form $T_{1}$ obtained in this way verifies the conditions from
the statement of the theorem.
$\qed$

As a corollary of the formula (\ref{T_2}) obtained in the proof above we can
make now the connection with the Lagrangian formalism. Namely, we have:

\begin{cor}

The expressions
$T_{\sigma}$
defined according to (\ref{T-sigma}) verify the generalized Helmholtz
equations. 
\label{helmholtz}
\end{cor}

{\bf Proof:}
We select from (\ref{Ts}) and (\ref{deltaT}) only those equations containing
the expressions
$T^{\emptyset}_{\sigma,i_{1},\dots,i_{n}}$.
They are 

\be
\partial^{I_{0}}_{\sigma_{0}} 
T^{\emptyset}_{\sigma_{1},i_{2},\dots,i_{n+1}} + 2 
(B_{0}T)^{I_{0},\emptyset}_{\sigma_{0},\sigma_{1},i_{2},\dots,i_{n+1}} = 0,
\quad |I_{0}| = s,
\label{H1}
\ee

\be
\partial^{I_{1}}_{\sigma_{1}} T^{\emptyset}_{\sigma_{2},i_{3},\dots,i_{n+2}} +
{1 \over n} \sum_{p=3}^{n+2} (-1)^{p-1} d_{i_{p}}
T^{I_{1},\emptyset}_{\sigma_{1},\sigma_{2},i_{3},\dots,
\hat{i_{p}},\dots,i_{n+2}} +
(B_{1}T)^{I_{1},\emptyset}_{\sigma_{1},\sigma_{2},i_{3},\dots,i_{n+2}},
\quad |I_{1}| = 1,\dots,s-1,
\label{H2}
\ee

\be
{1\over 2} \left(\partial_{\sigma_{1}} 
T^{\emptyset}_{\sigma_{2},i_{3},\dots,i_{n+2}} -
\partial_{\sigma_{2}} 
T^{\emptyset}_{\sigma_{1},i_{3},\dots,i_{n+2}} \right) +
{1 \over n} \sum_{p=3}^{n+2} (-1)^{p-1} d_{i_{p}}
T^{\emptyset,\emptyset}_{\sigma_{1},\sigma_{2},i_{3},\dots,
\hat{i_{p}},\dots,i_{n+2}} = 0.
\label{H3}
\ee

Next, one obtains from the formula (\ref{T_2}) and the antisymmetry property
(\ref{antisymmetry}) the following expression:

\be
T^{I_{1},\emptyset}_{\sigma_{1},\sigma_{2},i_{2},\dots,i_{n}} =
{1 \over 2} n \sum_{|J| \leq s-1-|I_{1}|-|I_{2}|} (-1)^{|J|+|I_{1}|}
{|J|+|I_{1}| \choose |J|} d_{J} \partial^{I_{1}Ji_{1}}_{\sigma_{2}}
T^{\emptyset}_{\sigma_{1},i_{1},\dots,i_{n}} + 
(\delta \Lambda)^{I_{1},\emptyset}_{\sigma_{1},\sigma_{2},i_{2},\dots,i_{n}}.
\label{T_2'}
\ee

Finally one substitutes the preceding formula into the equations (\ref{H1}) -
(\ref{H3}) and uses the well known identity:

\be
\sum_{p=0}^{n} (-1)^{p} \varepsilon_{i_{0},\dots,\hat{i}_{p},\dots,i_{n}} = 0.
\label{epsilon}
\ee

As a result one obtains the following set of equations

\be
\partial^{I}_{\sigma_{1}} T_{\sigma_{2}} = (-1)^{|I|} \sum_{|J| \leq s-|I|} 
(-1)^{|J|} {|J|+|I| \choose |J|} d_{J} 
\partial^{IJ}_{\sigma_{2}} T_{\sigma_{1}}
\label{HS}
\ee
for 
$|I| = s$,
$|I| = 1,\dots,s-1$
and
$|I| = 0$
respectively. But (\ref{HS}) are exactly the Helmholtz equations (see
\cite{AM}, \cite{AD}, \cite{AH}, \cite{GM}, \cite{Kr4}).
$\qed$

From the theorem above we also obtain:

\begin{prop}

The decomposition (\ref{decomp}) proved in the preceding theorem determines
in an unique way the form 
$T_{0}.$
\label{T=0}
\end{prop}

{\bf Proof:}
It is sufficient to prove that
$$
\alpha = T_{0} + dT_{1} = 0 \Longrightarrow T_{0} = 0.
$$

Indeed, from 
$\alpha = 0$
we have in particular
$T_{\sigma,i_{1},\dots,i_{n}} = 0$
and it follows from the formula (\ref{T-polyn}) that we have 
$T_{0} = 0.$
$\qed$

As a corollary, let us denote by
$[\alpha]$
the equivalence class of the form
$\alpha \in \Omega^{s,s, Lep}_{n+1,\xi,tr}(X)$
modulo
$
d\Omega^{s,s-2}_{n,\xi,2}(X) \cap \Omega^{s,s-1}_{n+1,\xi,2}(X).
$
Then we have

\be
[\alpha] = 0 \Longleftrightarrow T_{\sigma} = 0
\label{alpha-T}
\ee
which says that the class of $\alpha$ is uniquely determined by the so-called
{\it Euler-Lagrange components} of $\alpha$:
$T_{\sigma}, \sigma = 1, \dots, m \equiv N - n.$

We call the globally defined class 
$[\alpha]$
of a certain form
$\alpha \in \Omega^{s,s,Lep}_{n+1,\xi,tr}(X)$
a {\it Lagrange-Souriau class}. We note in closing this Section that there are
two particular but important cases when the class of the form $\alpha$ is
formed only from the form $\alpha$; obviously this happens when 

\be
d\Omega^{s,s-2}_{n,\xi,2}(X) \cap \Omega^{s,s-1}_{n+1,\xi,2}(X)
= 0.
\label{particular}
\ee

One can see that if 
$s = 2$, and $n$ 
arbitrary, or 
$n = 1$, and $s$ 
arbitrary, the equality above becomes an identity. So in this cases one can
speak of a globally defined {\it Lagrange-Souriau form} as in \cite{Gr1},
\cite{GP} and \cite{Gr2} respectively. 

\newpage
\section{The Associated Poincar\'e-Cartan $n$-Form}

\subsection{The General Construction}

Because, according to Corollary \ref{lin}, the Euler-Lagrange expressions 
$T_{\nu}$
are at most linear in the higher-order derivatives
$y^{\sigma}_{i_{1},\dots,i_{s}}$
it is to be expected that they follows from a Lagrangian of {\it minimal order}
$r \equiv {\left[ s + 1\right] \over 2}.$
We prove this fact in this Section. The key observation is

\begin{prop}

Let 
$\alpha \in \Omega^{s,s}_{q+1,\xi}(X)$
be closed. Then one can write it {\bf locally} in the form

\be
\alpha = d\theta
\label{a-t}
\ee
where
$\theta \in \Omega^{s}_{n}(X)$
is a $\rho^{s,s-1}_{n}$-projectable form and has the coordinates expression:

\be
\theta = \sum_{k=0}^{q} \sum_{|I_{1}|,\dots,|I_{k}| \leq r-1} 
L^{I_{1},\dots,I_{k}}_{\sigma_{1},\dots,\sigma_{k},i_{k+1},\dots,i_{q}}
\omega^{\sigma_{1}}_{I_{1}} \wedge \dots \wedge 
\omega^{\sigma_{k}}_{I_{k}} \wedge dx^{i_{k+1}} \wedge \dots \wedge dx^{i_{q}}
\label{theta}
\ee
where
$
L^{I_{1},\dots,I_{k}}_{\sigma_{1},\dots,\sigma_{k},i_{k+1},\dots,i_{q}},
\quad |I_{1}|,\dots,|I_{k}| \leq r - 1
$
are smooth functions depending on the variables
$(x^{i},y^{\sigma},y^{\sigma}_{j},\dots,y^{\sigma}_{j_{1},\dots,j_{s}})$
and verify the (anti)symmetry property (\ref{antisymmetry}).
\label{alpha-dtheta}
\end{prop}

{\bf Proof:}
Let us use the proposition \ref{project} and write 
$\alpha \in \Omega^{s,s}_{q+1,\xi}(X)$
in the form (\ref{alpha0}). It is easy to see that 
$\alpha_{0}$
has the following generic form:
$$
\alpha_{0} = \sum_{k=0}^{q+1} \sum_{|I_{1}|+\dots+|I_{k}| \leq s-1} 
A^{I_{1},\dots,I_{k}}_{\sigma_{1},\dots,\sigma_{k},i_{k+1},\dots,i_{q+1}}
dy^{\sigma_{1}}_{I_{1}} \wedge \dots \wedge dy^{\sigma_{k}}_{I_{k}} 
\wedge dx^{i_{k+1}} \wedge \dots \wedge dx^{i_{q+1}}
$$
where 
$
A^{I_{1},\dots,I_{k}}_{\sigma_{1},\dots,\sigma_{k},i_{k+1},\dots,i_{q+1}},
\quad |I_{1}|,\dots,|I_{k}| \leq s - 1
$
are smooth functions depending on the variables
$(x^{i},y^{\sigma},y^{\sigma}_{j},\dots,y^{\sigma}_{j_{1},\dots,j_{s}})$
and verify the (anti)symmetry property (\ref{antisymmetry}). In particular the
differentials 
$dy^{\sigma}_{I}, \quad |I| \geq r$
can appear at most once in every term of the preceding sum because if there 
would exists a term with at least two such differentials we would get a
contradiction according to the obvious inequality 
$2r \geq s.$
So, one can write

\be
\alpha_{0} = \beta + \sum_{|I| \geq r} dy^{\sigma}_{I} \wedge
\alpha^{I}_{\sigma}
\label{alpha-r}
\ee
where
$\beta$ and $\alpha^{I}_{\sigma}$ are $(q+1)$ (resp. $q$) forms which do not
contain the differentials
$dy^{\sigma}_{I}, \quad |I| \geq r.$

But the form $\alpha_{0}$ is closed (see prop. \ref{project}) so one can write
it, locally, as follows:

\be
\alpha_{0} = d\theta_{0}
\label{alpha-theta}
\ee
with $\theta_{0}$ having a structure similar to (\ref{alpha-r}):

\be
\theta_{0} = \gamma + \sum_{|I| \geq r} dy^{\sigma}_{I} \wedge
\theta^{I}_{\sigma};
\label{theta-r}
\ee
here
$\gamma$ and $\theta^{I}_{\sigma}$ are $q$ (resp. $(q-1)$) forms which do not
contain the differentials
$dy^{\sigma}_{I}, \quad |I| \geq r.$
If we substitute (\ref{alpha-r}) and (\ref{theta-r}) into (\ref{alpha-theta})
we easily obtain the consistency condition:
$$
\partial^{J}_{\nu} \theta^{I}_{\sigma} = 
\partial^{I}_{\sigma} \theta^{J}_{\nu}, \quad \forall |I|, |J| \geq r.
$$

Applying the usual Poincar\'e lemma one gets from here that
$\theta^{I}_{\sigma}, \quad |I| \geq r$
have the following expression
$$
\theta^{I}_{\sigma} = \partial^{I}_{\sigma} \lambda, \quad |I| \geq r
$$
where $\lambda$ is a $(q-1)$-form which do not
contain the differentials
$dy^{\sigma}_{I}, \quad |I| \geq r.$
Now one substitute this into the expression (\ref{theta-r}) above and gets:
$$
\theta_{0} = \gamma + d\lambda - \sum_{|I| \leq r} dy^{\sigma}_{I} \wedge
\partial^{I}_{\sigma} \lambda - dx^{i} \wedge {\partial \lambda \over \partial
x^{i}}.
$$

It follows that one can take in (\ref{alpha-theta}) 
$$
\theta_{0} = \gamma  - \sum_{|I| \leq r} dy^{\sigma}_{I} \wedge
\partial^{I}_{\sigma} \lambda - dx^{i} \wedge {\partial \lambda \over \partial
x^{i}}
$$
without affecting it. But it is clear that this form has the structure
$$
\theta_{0} = \sum_{k=0}^{q} \sum_{|I_{1}|+\dots+|I_{k}| \leq r-1} 
B^{I_{1},\dots,I_{k}}_{\sigma_{1},\dots,\sigma_{k},i_{k+1},\dots,i_{q}}
dy^{\sigma_{1}}_{I_{1}} \wedge \dots \wedge dy^{\sigma_{k}}_{I_{k}} 
\wedge dx^{i_{k+1}} \wedge \dots \wedge dx^{i_{q}}.
$$

If we define
$\theta \equiv (\rho^{s,s-1}_{n})^{*} \theta_{0}$
then we have the equality from the statement.
$\qed$

If $q = n$ then we have similarly to (\ref{T-sigma}):

\be
L_{i_{1},\dots,i_{n}} = \varepsilon_{i_{1},\dots,i_{n}} L
\label{lagr}
\ee
where $L$ is a smooth real  function called the {\it (local) Lagrangian}.

Applying prop \ref{alpha-dalpha} one gets now

\begin{cor}
Let $\alpha$ be as above. If we write it in the form (\ref{a-t}) with $\theta$
given by (\ref{theta}), then the following relations are true:

\be
\partial^{I_{0}}_{\sigma_{0}} 
L^{I_{1},\dots,I_{k}}_{\sigma_{1},\dots,\sigma_{k},i_{k+1},\dots,i_{q}} = 0
\quad |I_{0}| = s, \quad |I_{1}|,\dots,|I_{k}| \leq r-1, 
\quad k = 0,\dots,q,
\label{Ts1}
\ee

\begin{eqnarray}
\partial^{I_{0}}_{\sigma_{0}} 
L^{I_{1},\dots,I_{k}}_{\sigma_{1},\dots,\sigma_{k},i_{k+1},\dots,i_{q}} = 0
\quad |I_{0}| = r+1,\dots,s-1, \quad |I_{1}|,\dots,|I_{k}| \leq r-1, 
\nonumber \\
\quad |I_{0}|+ \dots + |I_{k}| \geq s, \quad k = 1,\dots,q,
\label{Ts2}
\end{eqnarray}

\begin{eqnarray}
\partial^{I_{0}}_{\sigma_{0}} 
L^{I_{1},\dots,I_{k}}_{\sigma_{1},\dots,\sigma_{k},i_{k+1},\dots,i_{q}} +(k+1) 
(B_{0}L)^{I_{0},\dots,I_{k}}_{\sigma_{0},\dots,\sigma_{k},i_{k+1},\dots,i_{q}}
= 0, \quad |I_{0}| = r,
\nonumber \\  \quad |I_{1}|,\dots,|I_{k}| \leq r-1, 
\quad |I_{0}|+ \dots + |I_{k}| \geq s, \quad k = 1,\dots,q,
\label{Ts3}
\end{eqnarray}

\be
(\delta L)^{I_{1},\dots,I_{k}}_{\sigma_{1},\dots,\sigma_{k},
i_{k+1},\dots,i_{q+1}} = 0, \quad |I_{1}|,\dots,|I_{k}| \leq s-1, 
\quad |I_{1}|+ \dots + |I_{k}| \geq s, \quad k = 2,\dots,q+1.
\label{deltaT1}
\ee

Moreover one has:

\be
\alpha = \sum_{k=0}^{q+1} \sum_{|I_{1}|,\dots,|I_{k}| \leq r} 
(\delta L)^{I_{1},\dots,I_{k}}_{\sigma_{1},\dots,\sigma_{k},
i_{k+1},\dots,i_{q+1}}
\omega^{\sigma_{1}}_{I_{1}} \wedge \dots \wedge \omega^{\sigma_{k}}_{I_{k}} 
\wedge dx^{i_{k+1}} \wedge \dots \wedge dx^{i_{q+1}}.
\label{a}
\ee
\label{cond-theta}
\end{cor}

Now we make use of the Lepage condition (\ref{T2}).

\begin{prop}

Let
$\alpha \in \Omega^{s,s}_{q+1,\xi}(X)$
be closed and verifying the Lepage condition (\ref{T2}). Suppose that we have
written it as in the proposition above. Then we have:

\be
\partial^{I}_{\sigma} L_{i_{1},\dots,i_{q}} = 0,\quad |I| > r.
\ee
\end{prop}

{\bf Proof:}
According to the preceding corollary we have: 

\be
T^{I}_{\sigma,i_{1},\dots,i_{q}} = (\delta L)^{I}_{\sigma,i_{1},\dots,i_{q}} = 
\partial^{I}_{\sigma} L_{i_{1},\dots,i_{q}} + {1 \over q} \sum_{p=1}^{q}
(-1)^{p} d_{i_{p}} L^{I}_{\sigma,i_{1},\dots,\hat{i_{p}},\dots,i_{q}} + 
(BL)^{I}_{\sigma,i_{1},\dots,i_{q}}, \quad |I| \leq s-1.
\label{lepage}
\ee

We take here 
$|I| > r,$
use the Lepage condition and also the fact that the the proposition
\ref{alpha-dtheta} implies:
$$
L^{I_{1}}_{\sigma_{1},i_{2},\dots,i_{q}} = 0, \quad |I_{1}| > r.
$$
The relation from the statement follows.
$\qed$

If $q = n$ then the same assertion is obviously true for the local Lagrangian 
$L$.

The Lepage condition is not completely exploited and can be used to provide an
expression for the local functions
$
L^{I_{1}}_{\sigma_{1},i_{2},\dots,i_{n}}, \quad |I| \leq  r - 1.
$

\begin{prop}

Let 
$\alpha \in \Omega^{s,s}_{n+1,\xi}(X)$
be closed and verifying the Lepage condition. Let us suppose that one writes 
$\alpha$ according to the proposition \ref{alpha-dtheta}. Then one has:

\be
L^{I_{1}}_{\sigma_{1},i_{2},\dots,i_{n}} = n \sum_{|J| \leq r-1-|I_{1}|}
(-1)^{|J|} d_{J} \partial^{I_{1}Ji_{1}}_{\sigma_{1}} L_{i_{1},\dots,i_{n}} +
(\delta \Lambda)^{I_{1}}_{\sigma_{1},i_{2},\dots,i_{n}}, 
\quad |I_{1}| \leq  r - 1
\label{L-I-sigma}
\ee
where we assume that the tensors 
$\Lambda^{\dots}_{\dots}$
have appropriate antisymmetry properties. We also admit that:

\be
\Lambda_{i_{2},\dots,i_{n}} = 0, 
\quad
\Lambda^{I_{1}}_{\sigma_{1},i_{2},\dots,i_{n-1}} = 0, 
\quad |I_{1}| \geq  r - 1,
\ee

and

\be
\partial^{I_{0}}_{\sigma_{0}} \Lambda^{I_{1}}_{\sigma_{1},i_{2},\dots,i_{n-1}} 
= 0, \quad |I_{0}| \geq 2r-1-|I_{1}|, \quad |I_{1}| \geq  r-2.
\ee
\label{L1}
\end{prop}

{\bf Proof:}
We take in the relation (\ref{lepage}) 
$1 \leq |I| \leq r$
and use the Lepage condition (\ref{T2}). The following set of equations
follows: 

\be
(BL)^{I}_{\sigma,i_{1},\dots,i_{n}} + {1 \over n} \sum_{p=1}^{n}
(-1)^{p} d_{i_{p}} L^{I}_{\sigma,i_{1},\dots,\hat{i_{p}},\dots,i_{n}} = 
- \partial^{I}_{\sigma} L_{i_{1},\dots,i_{n}}, \quad 1 \leq |I| \leq r.
\ee

From this equations one can determine recurrsively the expressions
$
L^{I_{1}}_{\sigma_{1},i_{2},\dots,i_{n}} 
$
for
$|I| = r - 1, r - 2, \dots, 0.$

We consider the system above as an inhomogeneous linear system of equations
(with the inhomogeneity determined in terms of the local Lagrangian) and
determine its solution as follows. First we prove that the first term in the
formula (\ref{L-I-sigma}) is a particular solution of this equation; this 
amounts to a simple computation using formula (\ref{epsilon}). Next, we
determine the general solution of {\it homogeneous} system of equation:

\be
(BL)^{I}_{\sigma,i_{1},\dots,i_{n}} + {1 \over n} \sum_{p=1}^{n}
(-1)^{p} d_{i_{p}} 
L^{I}_{\sigma,i_{1},\dots,\hat{i_{p}},\dots,i_{n}} = 0,
\quad 1 \leq |I| \leq r.
\ee

It is not very hard to prove by recurrence that the solution of this system of
equations is given by the second term in (\ref{L-I-sigma}).
$\qed$

Now we have a result similar to theorem 2 from \cite{Kr9}, ch. 3.2, but as we
can see the proof is much simpler and do not make use of the Young diagrams
technique. 

\begin{cor}

Let
$\alpha \in \Omega^{s,s}_{n+1,\xi}(X)$
be closed and verifying the Lepage condition (\ref{T2}). Suppose that we have
written it as in prop. \ref{alpha-dtheta}. Then the following formula is true:

\be
\theta = \theta_{0} + d\lambda + \mu
\ee
where

\be
\theta_{0} \equiv n! L dx^{1} \wedge \cdots \wedge dx^{n} +
n \sum_{|J| \leq r-1-|I_{1}|}
(-1)^{|J|} d_{J} \partial^{I_{1}Ji_{1}}_{\sigma_{1}} L_{i_{1},\dots,i_{n}} 
\omega^{\sigma_{1}}_{I_{1}} \wedge dx^{i_{2}} \wedge \cdots \wedge dx^{i_{n}};
\ee
also
$\lambda \in \Omega^{s,r-2}_{n-1,\xi,1}$
and
$\mu \in \Omega^{s,s-1}_{n,\xi,2}.$
\end{cor}

{\bf Proof:}
We start from the formula (\ref{theta}) obtained before and notice that the 
term corresponding to $k = 1$ can be written, according to the propositions
\ref{alpha-dalpha} and \ref{L1} as follows:
$$
L^{I_{1}}_{\sigma_{1},i_{2},\dots,i_{n}} 
\omega^{\sigma_{1}}_{I_{1}} \wedge dx^{i_{2}} \wedge \cdots \wedge dx^{i_{n}} 
= n \sum_{|J| \leq r-1-|I_{1}|} (-1)^{|J|} d_{J} 
\partial^{I_{1}Ji_{1}}_{\sigma_{1}} L_{i_{1},\dots,i_{n}} 
\omega^{\sigma_{1}}_{I_{1}} \wedge dx^{i_{2}} \wedge \cdots \wedge dx^{i_{n}} +
d\lambda
$$
where
$$
\lambda \equiv \sum_{|I_{1}| \leq r-2}
\Lambda^{I_{1}}_{\sigma_{1},i_{2},\dots,i_{n-1}}
\omega^{\sigma_{1}}_{I_{1}} \wedge dx^{i_{2}} \wedge \cdots 
\wedge dx^{i_{n-1}}. 
$$
is a $(n-1)$-form with the order of contactness equal to $1$.
Now we define the form $\mu$ to be the sum of the terms corresponding to the
contributions $k  \geq 2$ in the expression (\ref{theta}); this gives us a
$n$-form with the order of contactness equal to $2$. The formula from
the statement follows.
$\qed$

We prove now that the Euler-Lagrange expressions 
$T_{\sigma}$
are following from a Lagrangian of order $r$. 

\begin{prop}

In the conditions of prop. \ref{L1} the following result is true:

\be
T_{\sigma} = \sum_{|J| \leq r} (-1)^{|J|} d_{J} \partial^{J}_{\sigma} L.
\ee
\end{prop}

{\bf Proof:}
We have from (\ref{lepage}) for $I = \emptyset$ and $q = n$:

\be
T^{\emptyset}_{\sigma,i_{1},\dots,i_{n}} = 
\partial_{\sigma} L_{i_{1},\dots,i_{n}} + {1 \over n} \sum_{p=1}^{n}
(-1)^{p} d_{i_{p}} L^{\emptyset}_{\sigma,i_{1},\dots,\hat{i_{p}},\dots,i_{n}}. 
\ee

We also have from proposition \ref{L1} for $I_{1} = \emptyset$:

\be
L^{\emptyset}_{\sigma_{1},i_{2},\dots,i_{n}} = n \sum_{|J| \leq r-1}
(-1)^{|J|} d_{J} \partial^{Ji_{1}}_{\sigma_{1}} L_{i_{1},\dots,i_{n}} +
(\delta \Lambda)^{\emptyset}_{\sigma_{1},i_{2},\dots,i_{n}}.
\ee

If we substitute the second relation into the first one, we obtain the formula
from the statement.
$\qed$

Let us comment this result. First we can say that because the expressions
$T_{\sigma}$ 
have the usual Euler-Lagrange expression, they verify the generalized Helmholtz
equations, so we have an alternative proof of the corollary \ref{helmholtz}.
Next, we notice that in fact we have a sharper result, namely the expressions
$T_{\sigma}$ 
follow from a Lagrangian of order $r$ which is the minimal possible order.
Indeed, if 
$T_{\sigma}$ 
would follow from a Lagrangian of order strictly smaller than $r$, the the
Euler-Lagrange equations would have the order strictly smaller than $s$ which
would contradict the basic stating point of our analysis. So, we can say that
we have obtained above a form of the conjecture regarding the {\it reduction 
to the minimal order} in the higher-order Lagrangian formalism. 

We close this Subsection with the following result.

\begin{prop}

In the conditions of the proposition \ref{alpha-dtheta}, let us suppose that
$q = n$ and moreover that the tensors
$L^{I_{1},\dots,I_{k}}_{\sigma_{1},\dots,\sigma_{k},i_{k+1},\dots,i_{n}},
\quad |I_{1}| = \cdots = |I_{k}| = r-1, \quad k =1,\dots,n$
are traceless. The we have the following formula:

\begin{eqnarray}
L^{I_{1},\dots,I_{k}}_{\sigma_{1},\dots,\sigma_{k},i_{k+1},\dots,i_{n}} =
{n \choose k} r^{k} {(r-1)! \over (k+r-1)!} 
\partial^{i_{1}I_{1}}_{\sigma_{1}} \cdots \partial^{i_{k}I_{k}}_{\sigma_{k}} 
L_{i_{1},\dots,i_{n}}, 
\nonumber \\ \quad |I_{1}| = \cdots = |I_{k}| = r-1, \quad k =0,\dots,n
\end{eqnarray}
\end{prop}

The proof goes by induction and is based on the relation (\ref{Ts3}) and the
condition (\ref{T3}). We do not give the details but we only mention that the
preceding formula appears also in \cite{Kr5}.

\subsection{Hyper-Jacobians and Variationally Trivial Lagrangians}

We have emphasized before that equations of the type (\ref{Ts}) (see also
(\ref{Ts3})) play an important r\^ole in the Lagrangian formalism. Mainly, they
provide a natural way of obtaining some polynomial objects called
hyper-Jacobians \cite{BCO}, \cite{Ol1}. Let us prove this fact. By definition,
the {\it hyper-Jacobians of order s} are the following expressions:

\be
{\cal J}^{\sigma_{1},\dots,\sigma_{k},i_{k+1},\dots,i_{n}}_{I_{1},\dots,I_{k}} 
\equiv \varepsilon^{i_{1},\dots,i_{n}} \prod_{l=1}^{k}
y^{\sigma_{l}}_{I_{l}i_{l}}, \quad |I_{1}| = \cdots = |I_{k}| = s-1, 
\quad k = 0, \dots, n.
\label{jacobians}
\ee 

Then we have a series of results of combinatorial nature.

\begin{thm}

Let us consider the tensors
$
L_{k} \equiv \left( 
L^{I_{1},\dots,I_{k}}_{\sigma_{1},\dots,\sigma_{k},i_{k+1},\dots,i_{q}}\right),
\quad |I_{1}|,\dots,|I_{k}| \leq s-1
$
depending on the variables
$(x^{i},y^{\sigma},y^{\sigma}_{j},\dots,y^{\sigma}_{j_{1},\dots,j_{s}})$
and verifying the symmetry properties (\ref{antisymmetry}). Then these tensors
verify the system of equations:

\begin{eqnarray}
\partial^{I_{0}}_{\sigma_{0}} 
L^{I_{1},\dots,I_{k}}_{\sigma_{1},\dots,\sigma_{k},i_{k+1},\dots,i_{q}} +(k+1) 
(B_{0}L)^{I_{0},\dots,I_{k}}_{\sigma_{0},\dots,\sigma_{k},i_{k+1},\dots,i_{q}}
= 0, \nonumber \\ \quad |I_{0}| = s, \quad |I_{1}|,\dots,|I_{k}| \leq s-1, 
\quad k = 0,\dots,q
\label{Ls}
\end{eqnarray}

if and only if they have the following form:

-for $q \leq n$ and $k = 0,\dots, q:$

\begin{eqnarray}
L^{I_{1},\dots,I_{k}}_{\sigma_{1},\dots,\sigma_{k},i_{k+1},\dots,i_{q}} =
(-1)^{k(q+1)} {1 \over n!} {n \choose q-k} \varepsilon_{i_{1},\dots,i_{n}}
\sum_{l=k}^{q} {l \choose k} 
\nonumber \\
\sum_{|I_{k+1}|=\cdots =|I_{l}|=s-1}
{\cal L}^{I_{1},\dots,I_{l}}_{\sigma_{1},\dots,\sigma_{l},j_{l+1},\dots,j_{q}}
{\cal J}^{\sigma_{k+1},\dots,\sigma_{l},j_{l+1},\dots,j_{q},i_{1},\dots,i_{k},
i_{q+1},\dots,i_{n}}_{I_{k+1},\dots,I_{l}}, 
\label{L-jacob-q}
\end{eqnarray}

- for $q \geq n$ and $k = q-n,\dots,q:$

\begin{eqnarray}
L^{I_{1},\dots,I_{k}}_{\sigma_{1},\dots,\sigma_{k},i_{k+1},\dots,i_{q}} =
(-1)^{k(n+1)} {1 \over n!} {n \choose q-k} \varepsilon_{i_{q-n+1},\dots,i_{q}}
\sum_{l=k}^{q} {l \choose k} 
\nonumber \\ 
\sum_{|I_{k+1}|=\cdots =|I_{l}|=s-1}
{\cal L}^{I_{1},\dots,I_{l}}_{\sigma_{1},\dots,\sigma_{l},j_{l+1},\dots,j_{q}}
{\cal J}^{\sigma_{k+1},\dots,\sigma_{l},j_{l+1},\dots,j_{q},
i_{q-n+1},\dots,i_{k}}_{I_{k+1},\dots,I_{l}}. \quad
\label{L-jacob-n}
\end{eqnarray}

Here 
$
{\cal L}^{I_{1},\dots,I_{l}}_{\sigma_{1},\dots,\sigma_{l},j_{l+1},\dots,j_{q}}
$
are tensors depending on the variables
$(x^{i},y^{\sigma},y^{\sigma}_{j},\dots,y^{\sigma}_{j_{1},\dots,j_{s-1}})$
with the symmetry properties (\ref{antisymmetry}) and
$
{\cal J}^{\sigma_{1},\dots,\sigma_{l},j_{l+1},\dots,j_{q}}_{I_{1},\dots,I_{l}}
$
are the hyper-Jacobians of order $s$.
\label{L-jacob}
\end{thm}

{\bf Proof:}
The proof is based on some tedious computations, but it is straightforward so
we omit many details. We present only the case $q \leq n$, the other one being
dealt with similarly.  Let us start from the equation in the statement
(\ref{Ls}) for $k = q$:
$$
\partial^{I_{0}}_{\sigma_{0}} 
L^{I_{1},\dots,I_{q}}_{\sigma_{1},\dots,\sigma_{q}} = 0
$$
i.e. 
$L^{I_{1},\dots,I_{q}}_{\sigma_{1},\dots,\sigma_{q}}$
does not depend on 
$y^{\sigma}_{j_{1},\dots,j_{s}}.$
So, the formula (\ref{L-jacob-q}) is true for $k = q$ i.e. we have
$$
L^{I_{1},\dots,I_{q}}_{\sigma_{1},\dots,\sigma_{q}} = - {1 \over n!} 
\varepsilon_{i_{1},\dots,i_{n}} 
{\cal L}^{I_{1},\dots,I_{q}}_{\sigma_{1},\dots,\sigma_{q}} 
{\cal J}^{i_{1},\dots,i_{n}} .
$$

Now we proceed by induction downward on $k$. We suppose that the formula
(\ref{L-jacob-q}) is true for 
$k+1,\dots,q$
and we prove it for $k.$ One considers the equations (\ref{Ls}) and notices
that the right hand side is known from the induction hypothesis. So, this
system can be considered as an inhomogeneous system for the unknown variables
$L^{I_{1},\dots,I_{k}}_{\sigma_{1},\dots,\sigma_{k},i_{k+1},\dots,i_{q}}.$
After some computations this system proves to be:
\begin{eqnarray}
\partial^{I_{0}}_{\sigma_{0}} 
L^{I_{1},\dots,I_{k}}_{\sigma_{1},\dots,\sigma_{k},i_{k+1},\dots,i_{q}} =
- (-1)^{k(q+1)} (k+1) {1 \over n!} {n \choose q-k} 
\varepsilon_{i_{1},\dots,i_{n}} \sum_{l=k}^{q-1} {l+1 \choose k+1} 
\nonumber \\
\sum_{|I_{k+1}|=\cdots |I_{l}|=s-1}
(B_{0}{\cal L})^{I_{0},\dots,I_{l}}_{\sigma_{0},\dots,\sigma_{l},
j_{l+1},\dots,j_{q}} 
\times {\cal J}^{\sigma_{k+1},\dots,\sigma_{l},j_{l+1},\dots,j_{q},i_{1},\dots,
i_{k},i_{q+1},\dots,i_{n}}_{I_{k+1},\dots,I_{l}}. 
\end{eqnarray}

The general solution of this system is a sum of the general solution of the
homogeneous system:
$$
\partial^{I_{0}}_{\sigma_{0}} 
L^{I_{1},\dots,I_{k}}_{\sigma_{1},\dots,\sigma_{k},i_{k+1},\dots,i_{q}} = 0
$$
(i.e. a tensor independent of
$y^{\sigma}_{j_{1},\dots,j_{s}}$)
and a particular solution of the inhomogeneous system. So it remains to find
out such a particular solution. As expected, such a particular solution is
given exactly by the expression for
$
L^{I_{1},\dots,I_{k}}_{\sigma_{1},\dots,\sigma_{k},i_{k+1},\dots,i_{q}}
$
from the statement. The induction is finished and we have proven that the
system (\ref{Ls}) has the general solution given by (\ref{L-jacob-q}). From the
proof it is clear that the converse is also true. The case $q \geq n$ is
considered along similar lines.
$\qed$

We have succeeded in this way to obtain the hyper-Jacobians as solution of a
certain system of equations derived from a closedness condition (see cor.
\ref{lin}). The following consequence follows.

\begin{cor}

In the conditions of proposition \ref{alpha-dtheta} we have for 
$k \geq 2$
and
$|I_{1}| + \dots + |I_{k}| \geq  r$
the following formulae:

-for $q \leq n$:

\begin{eqnarray}
L^{I_{1},\dots,I_{k}}_{\sigma_{1},\dots,\sigma_{k},i_{k+1},\dots,i_{q}} =
(-1)^{k(q+1)} {1 \over n!} {n \choose q-k} \varepsilon_{i_{1},\dots,i_{n}}
\sum_{l=k}^{q} {l \choose k} 
\nonumber \\
\sum_{|I_{k+1}|=\cdots =|I_{l}|=r-1}
{\cal L}^{I_{1},\dots,I_{l}}_{\sigma_{1},\dots,\sigma_{l},j_{l+1},\dots,j_{q}}
{\cal J}^{\sigma_{k+1},\dots,\sigma_{l},j_{l+1},\dots,j_{q},i_{1},\dots,i_{k},
i_{q+1},\dots,i_{n}}_{I_{k+1},\dots,I_{l}};
\end{eqnarray}

- for $q \geq n$:

\begin{eqnarray}
L^{I_{1},\dots,I_{k}}_{\sigma_{1},\dots,\sigma_{k},i_{k+1},\dots,i_{q}} =
(-1)^{k(n+1)} {1 \over n!} {n \choose q-k} \varepsilon_{i_{q-n+1},\dots,i_{q}}
\sum_{l=k}^{q} {l \choose k} 
\nonumber \\
\sum_{|I_{k+1}|=\cdots =|I_{l}|=r-1}
{\cal L}^{I_{1},\dots,I_{l}}_{\sigma_{1},\dots,\sigma_{l},j_{l+1},\dots,j_{q}}
{\cal J}^{\sigma_{k+1},\dots,\sigma_{l},j_{l+1},\dots,j_{q},
i_{q-n+1},\dots,i_{k}}_{I_{k+1},\dots,I_{l}} \quad
\end{eqnarray}
where 
${\cal J}^{\sigma_{k+1},\dots,\sigma_{l},
j_{l+1},\dots,j_{q}}_{I_{k+1},\dots,I_{l}}$
are the hyper-Jacobians of order $r$ and
$
{\cal L}^{I_{1},\dots,I_{l}}_{\sigma_{1},\dots,\sigma_{l},j_{l+1},\dots,j_{q}},
\quad |I_{1}|,\dots,|I_{k}| \leq r-1
$
are tensors depending on the variables
$(x^{1},y^{\sigma},\dots,y^{\sigma}_{j_{1},\dots,j_{r-1}})$
and with the symmetry properties (\ref{antisymmetry}).
If $s = 2r-1$ then the formulae above are valid for $k = 1$ also.
\end{cor}

{\bf Proof:}
We use the first three formulae from the corollary \ref{cond-theta} and obtain
in the conditions of the statement:

\be
\partial^{I_{0}}_{\sigma_{0}} 
L^{I_{1},\dots,I_{k}}_{\sigma_{1},\dots,\sigma_{k},i_{k+1},\dots,i_{q}} = 0,
\quad |I_{0}| \geq r+1,
\ee

\be
\partial^{I_{0}}_{\sigma_{0}} 
L^{I_{1},\dots,I_{k}}_{\sigma_{1},\dots,\sigma_{k},i_{k+1},\dots,i_{q}} +(k+1) 
(B_{0}T)^{I_{0},\dots,I_{k}}_{\sigma_{0},\dots,\sigma_{k},i_{k+1},\dots,i_{q}}
= 0, \quad |I_{0}| = r.
\ee

Now we can apply the same arguments as in the preceding theorem.
$\qed$

In the following we prove that the expressions just derived (\ref{L-jacob-q})
and (\ref{L-jacob-n}) behave naturally with respect to the cohomology operator
$\delta$ defined by (\ref{delta12}). First we define, similarly
(\ref{delta12}):

\be
\delta' \equiv \delta_{1} + \delta'_{2}
\label{delta'}
\ee
where
$\delta'_{2}$
has the same structure as
$\delta_{2}$
given by (\ref{delta2}) but with
$d_{i} \mapsto d'_{i} \equiv d_{i}^{s-1}.$
Now we have:

\begin{thm}

In the conditions of the theorem above the following formulae are true:

-for $q \leq n - 1$ and $k = 0,\dots,q+1:$

\begin{eqnarray}
(\delta L)^{I_{1},\dots,I_{k}}_{\sigma_{1},\dots,\sigma_{k},
i_{k+1},\dots,i_{q+1}} = (-1)^{kq)} {1 \over n!} {n \choose q-k+1} 
\varepsilon_{i_{1},\dots,i_{n}}
\sum_{l=k}^{q} {l \choose k} 
\nonumber \\
\sum_{|I_{k+1}|=\cdots =|I_{l}|=s-1}
(\delta' {\cal L})^{I_{1},\dots,I_{l}}_{\sigma_{1},\dots,\sigma_{l},
j_{l+1},\dots,j_{q+1}} 
{\cal J}^{\sigma_{k+1},\dots,\sigma_{l},j_{l+1},\dots,j_{q+1},
i_{1},\dots,i_{k},i_{q+2},\dots,i_{n}}_{I_{k+1},\dots,I_{l}}; \quad
\label{delta-L-jacob-q} 
\end{eqnarray}

- for $q \geq n - 1$ and $k = q-n+1,\dots,q:$

\begin{eqnarray}
(\delta L)^{I_{1},\dots,I_{k}}_{\sigma_{1},\dots,\sigma_{k},
i_{k+1},\dots,i_{q+1}} = (-1)^{k(n-1)} {1 \over n!} {n \choose q-k+1} 
\varepsilon_{i_{q-n+2},\dots,i_{q+1}}
\sum_{l=k}^{q} {l \choose k} 
\nonumber \\
\sum_{|I_{k+1}|=\cdots |I_{l}|=s-1}
(\delta' {\cal L})^{I_{1},\dots,I_{l}}_{\sigma_{1},\dots,\sigma_{l},
j_{l+1},\dots,j_{q+1}}  
{\cal J}^{\sigma_{k+1},\dots,\sigma_{l},j_{l+1},\dots,j_{q+1},
i_{q-n+2},\dots,i_{k}}_{I_{k+1},\dots,I_{l}}. \quad
\label{delta-L-jacob-n} 
\end{eqnarray}
\label{natural}
\end{thm}

The proof consists from hard computations and the definitions of the various
expressions on the left hand side making use of the preceding theorem. Because
these computations are only a question of combinatorial ability, we skip them.

Based on the preceding two results we can give now a new proof of an important
result from \cite{Gr4}.

\begin{thm}

A (local) Lagrangian $L$ of order $r$ is variationally trivial (i.e. the
associated Euler-Lagrange expressions are identically zero) if and only if it
has the following form:

\be
L = \sum_{k=0}^{n} \sum_{|I_{1}|=\cdots =|I_{k}|=r-1}
(\delta' \lambda)^{I_{1},\dots,I_{k}}_{\sigma_{1},\dots,\sigma_{k},
i_{k+1},\dots,i_{n}} \times 
{\cal J}^{\sigma_{1},\dots,\sigma_{k},i_{k+1},\dots,i_{n}}_{I_{1},\dots,I_{k}}
\label{lagr-triv}
\ee
where 
$
\lambda^{I_{1},\dots,I_{k}}_{\sigma_{1},\dots,\sigma_{k},i_{k+1},\dots,i_{n-1}}
, \quad k = 0,\dots,n-1
$
are arbitrary tensors depending on the variables
$(x^{1},y^{\sigma},\dots,y^{\sigma}_{j_{1},\dots,j_{r-1}})$
and verifying the symmetry properties (\ref{antisymmetry}). Here
$
{\cal J}^{\sigma_{1},\dots,\sigma_{k},i_{k+1},\dots,i_{n}}_{I_{1},\dots,I_{k}}
$
are the hyper-Jacobians of order $r$.
\end{thm}

{\bf Proof:}
(i) According to the theorem  \ref{canonical} (more directly from the
proposition \ref{T=0}) we have that 
$T_{\sigma} = 0$
{\it iff}
$\alpha \in d \Omega^{s,s-2}_{n,\xi,2}(X) \cap \Omega^{s,s-1}_{n+1,\xi,2}(X)$.

More explicitly we have:

\be
\alpha = \sum_{k=2}^{n} \sum_{|I_{1}|+\cdots +|I_{k}| \leq s-2}
(\delta \Lambda)^{I_{1},\dots,I_{k}}_{\sigma_{1},\dots,\sigma_{k},
i_{k+1},\dots,i_{n+1}} \omega^{\sigma_{1}}_{I_{1}} \wedge \dots \wedge
\omega^{\sigma_{k}}_{I_{k}} \wedge dx^{i_{k+1}} \wedge \cdots \wedge 
dx^{i_{n+1}}
\label{alpha-trivial}
\ee
such that the equations (\ref{lambda1}) and (\ref{lambda2}) are verified by the
tensors
$\Lambda^{\dots}_{\dots}.$

If we compare with the standard expression of $\alpha$ given by (\ref{alpha}) 
we get from here in particular: 

\be
T^{I_{1}}_{\sigma_{1},i_{2},\dots,i_{n+1}} = 0, \quad |I_{1}| \leq s-1
\ee
and
\be
T^{I_{1},\dots,I_{k}}_{\sigma_{1},\dots,\sigma_{k},i_{k+1},\dots,i_{n+1}} = 0,
\quad |I_{1}| + \cdots |I_{k}| \geq s-1, \quad k = 2,\dots,n.
\ee

(ii) Now we assume that we have written $\alpha$ as in Theorem
\ref{alpha-dtheta}.  We compare this expression with the two previous relations
and we obtain the following restrictions on the tensors
$L^{I_{1},\dots,I_{k}}_{\sigma_{1},\dots,\sigma_{k},i_{k+1},\dots,i_{n}}:$ 

In the following we will consider only the case $s > 2$, the case $s = 2$
affording a direct and simple analysis. We obtain: 

- from the first relation:

\be
\partial^{I_{0}}_{\sigma_{0}} L_{i_{1},\dots,i_{n}} = 0, \quad 
|I_{0}| \geq r+1,
\label{1-r+1}
\ee

\be
\partial^{I_{0}}_{\sigma_{0}} L_{i_{1},\dots,i_{n}} +
(B_{0}L)^{I_{0}}_{\sigma_{0},i_{1},\dots,i_{n}} = 0, \quad |I_{0}| = r
\label{1-r}
\ee
and

\be
(\delta L)^{I_{1}}_{\sigma_{1},i_{2},\dots,i_{n+1}} = 0, 
\quad |I_{1}| \leq r-1.
\label{1-r-1}
\ee

- from the second relation:

\be
\partial^{I_{0}}_{\sigma_{0}} 
L^{I_{1},\dots,I_{k}}_{\sigma_{1},\dots,\sigma_{k},i_{k+1},\dots,i_{n}} = 0,
\quad |I_{0}| \geq r+1, \quad |I_{1}|, \dots, |I_{k}| \geq r-2, \quad k \geq 1,
\label{2-r+1}
\ee

\be
\partial^{I_{0}}_{\sigma_{0}} 
L^{I_{1}}_{\sigma_{1},i_{2},\dots,i_{n}} + 2 
(B_{0}L)^{I_{0},I_{1}}_{\sigma_{0},\sigma_{1},i_{2},\dots,i_{n}} = 0,
\quad |I_{0}| = r, \quad |I_{1}| = r-1,
\label{2-r-r-1}
\ee

\be
\partial^{I_{0}}_{\sigma_{0}} 
L^{I_{1},\dots,I_{k}}_{\sigma_{1},\dots,\sigma_{k},i_{k+1},\dots,i_{n}} +(k+1) 
(B_{0}L)^{I_{0},\dots,I_{k}}_{\sigma_{0},\dots,\sigma_{k},i_{k+1},\dots,i_{n}}
= 0, \quad |I_{0}| = r, \quad |I_{1}|, \dots, |I_{k}| \geq r-2, \quad k \geq 2,
\label{2-r-r-2}
\ee
and

\be
(\delta L)^{I_{1},\dots,I_{k}}_{\sigma_{1},\dots,\sigma_{k},
i_{k+1},\dots,i_{n+1}}, \quad |I_{1}| = \cdots = |I_{k}| = r-1,
\quad k = 3,\dots n+1.
\label{2-3}
\ee

Finally, if we take into account the properties of the tensors 
$\Lambda^{\dots}_{\dots}$
appearing (\ref{alpha-trivial}), namely (\ref{lambda1}) and (\ref{lambda2}), 
we have also:

\be
\partial^{I_{0}}_{\sigma_{0}} 
L^{I_{1}}_{\sigma_{1},i_{2},\dots,i_{n}} + 2 
(B_{0}L)^{I_{0},I_{1}}_{\sigma_{0},\sigma_{1},i_{2},\dots,i_{n}} = 
\left( (B_{0} + B_{1})\Lambda\right)^{I_{0},I_{1}}_{\sigma_{0},\sigma_{1},
i_{2},\dots,i_{n}}, \quad |I_{0}| = r, \quad |I_{1}| = r-2.
\label{2-r-r-2(1)}
\ee

(iii) We have to exploit these relations to obtain a generic form for $L$.
First we note that from (\ref{1-r+1}) it follows that $L$ depends only on the
variables 
$(x^{i},y^{\sigma},\dots,y^{\sigma}_{j_{1},\dots,j_{r}}).$ 
From (\ref{2-r+1}) it follows that the same assertion is true for the tensors
$L^{I_{1},\dots,I_{k}}_{\sigma_{1},\dots,\sigma_{k},i_{k+1},\dots,i_{n}}$ 
for
$I_{1}|, \dots, |I_{k}| \geq r-2, \quad k \geq 1.$

Now we can use the relations (\ref{2-r-r-2}) in the same way as in Theorem 
\ref{L-jacob} to obtain that:

\begin{eqnarray}
L^{I_{1},\dots,I_{k}}_{\sigma_{1},\dots,\sigma_{k},i_{k+1},\dots,i_{n}} =
(-1)^{k(n+1)} {1 \over n!} {n \choose k} \varepsilon_{i_{1},\dots,i_{n}}
\sum_{l=k}^{n} {l \choose k} 
\nonumber \\
\sum_{|I_{k+1}|=\cdots =|I_{l}|=r-1}
{\cal L}^{I_{1},\dots,I_{l}}_{\sigma_{1},\dots,\sigma_{l},j_{l+1},\dots,j_{n}}
{\cal J}^{\sigma_{k+1},\dots,\sigma_{l},j_{l+1},\dots,j_{n},
i_{1},\dots,i_{k}}_{I_{k+1},\dots,I_{l}}, \quad k = 2,\dots, q;
\label{k>2}
\end{eqnarray}
here the multi-indices verify
$|I_{1}|, \dots, |I_{k}| \geq r-2, \quad k \geq 2$
but at most one of them has the length $r-2$. 
As in Theorem \ref{L-jacob}, the expressions
$
{\cal L}^{I_{1},\dots,I_{l}}_{\sigma_{1},\dots,\sigma_{l},j_{l+1},\dots,j_{n}},
$
defined for
$|I_{1}|, \dots, |I_{k}| \geq r-2, \quad k \geq 2$
but at most one of them has the length $r-2,$
are tensors depending on the variables
$(x^{i},y^{\sigma},y^{\sigma}_{j},\dots,y^{\sigma}_{j_{1},\dots,j_{r-1}})$
with the symmetry properties (\ref{antisymmetry});
$
{\cal J}^{\sigma_{1},\dots,\sigma_{l},j_{l+1},\dots,j_{q}}_{I_{1},\dots,I_{l}}
$
are the hyper-Jacobians of order $r$.
If we use (\ref{2-r-r-1}) we can extend this relation for $k = 1$ and
$|I_{1}| = r-1$;
if we use (\ref{1-r}) we can extended it for $k = 0$ also. Explicitly:

\be
L = {1 \over n!} \varepsilon_{i_{1},\dots,i_{n}}
\sum_{l=0}^{n} \sum_{|I_{1}|=\cdots =|I_{l}|=r-1}
{\cal L}^{I_{1},\dots,I_{l}}_{\sigma_{1},\dots,\sigma_{l},j_{l+1},\dots,j_{n}}
{\cal J}^{\sigma_{1},\dots,\sigma_{l},j_{l+1},\dots,j_{n},}_{I_{1},\dots,I_{l}}
\label{L-tr}
\ee
which is already a big step in proving the formula from the statement. 

Because the right hand side of (\ref{2-r-r-2(1)}) is not zero, we cannot extend
the relation (\ref{k>2}) for $k = 1$ and 
$|I_{1}| = r-2.$ 
However, one can prove that one can modify the tensors
$L^{I_{1},\dots,I_{k}}_{\sigma_{1},\dots,\sigma_{k},i_{k+1},\dots,i_{n}}$
in such a way that the form $\alpha$ is not changed, but one has for $k = 1$ 
and
$|I_{1}| \leq r-2:$

\be
L^{I_{1}}_{\sigma_{1},i_{2},\dots,i_{n}} = 
L^{(0)I_{1}}_{\sigma_{1},i_{2},\dots,i_{n}} +
(\delta K)^{I_{1}}_{\sigma_{1},i_{2},\dots,i_{n}};
\ee
here the tensors
$L^{(0)I_{1}}_{\sigma_{1},i_{2},\dots,i_{n}}$
are given by an expression appearing in the right hand side of (\ref{k>2}),
we admit that
$K_{i_{1},\dots,i_{n-1}} = 0$
and that
$K^{I_{1}}_{\sigma,i_{2},\dots,i_{n-1}}, \quad |I_{1}| \leq r-3$
are arbitrary tensors depending on the variables
$(x^{i},y^{\sigma},y^{\sigma}_{j},\dots,y^{\sigma}_{j_{1},\dots,j_{r}}).$
One can group the contribution of the last term in the form $\theta$ into an
exact differential so one can neglect it. It follows that one has (\ref{k>2})
for $k = 1$ 
and
$|I_{1}| \leq r-2.$

Moreover, one can fix the function
${\cal L}^{\dots}_{\dots}$
such that the equation (\ref{2-3}) stays true for $k =2.$ This gives easily:

\be
(\delta' {\cal L})^{I_{1},I_{2}}_{\sigma_{1},\sigma_{2},j_{3},\dots,j_{n+1}} = 
0, \quad |I_{1}| = |I_{2}| = r-1.
\label{2-3-cal(2)}
\ee

We also have:

\be
(\delta' {\cal L})^{I_{1}}_{\sigma_{1},j_{2},\dots,j_{n+1}} = 0, 
\quad \forall |I_{1}| \leq r-1.
\label{2-3-cal(1)}
\ee

Finally if we insert (\ref{k>2}) into (\ref{2-3}) we obtain, using Theorem
\ref{natural}:

\be
(\delta' {\cal L})^{I_{1},\dots,I_{l}}_{\sigma_{1},\dots,\sigma_{l},
j_{l+1},\dots,j_{n+1}} = 0, \quad |I_{1}| = \cdots + |I_{l}| = r-1,
\quad k = 3,\dots,n+1
\label{2-3-cal}
\ee

(iv) It remains to find out the most general solution of the system of
equations (\ref{2-3-cal(2)}), (\ref{2-3-cal(1)}) and (\ref{2-3-cal}). This is a
cohomological problem which can be solved by a descent procedure of the same
type as the one appearing in the BRST cohomology. 

We start from (\ref{2-3-cal}) for $l = n+1:$
$$
(\delta' {\cal L})^{I_{1},\dots,I_{n+1}}_{\sigma_{1},\dots,\sigma_{n+1}} = 0
\quad \Leftrightarrow \quad
(\delta_{1} {\cal L})^{I_{1},\dots,I_{n+1}}_{\sigma_{1},\dots,\sigma_{n+1}} = 0
, \quad |I_{1}| = \dots = |I_{n+1}| = r-1
$$
or, explicitly:

\be
\sum_{p=1}^{n+1} (-1)^{p-1} \partial^{I_{p}}_{\sigma_{p}}
{\cal L}^{I_{1},\dots,\hat{I_{p}},\dots,I_{n+1}}_{\sigma_{1},\dots,
\hat{\sigma_{p}},\dots,\sigma_{n+1}} = 0, 
\quad |I_{1}| = \dots = |I_{n+1}| = r-1.
\ee

This type of relation can be easily solved using a de Rham type homology
operator: we define

\be
\lambda^{I_{1},\dots,I_{n-1}}_{\sigma_{1},\dots,\sigma_{n-1}} \equiv
\int_{0}^{1} t^{n-1} \sum_{I_{0}|=r-1} y^{\sigma_{0}}_{I_{0}}
{\cal L}^{I_{0},\dots,I_{n-1}}_{\sigma_{0},\dots,\sigma_{n-1}}
(x^{y},y^{\sigma},\dots,ty^{j_{1},\dots,j_{r-1}}) dt 
\ee
and obtain that

\be
{\cal L}^{I_{1},\dots,I_{n}}_{\sigma_{1},\dots,\sigma_{n}} =
(\delta_{1} {\cal \lambda})^{I_{1},\dots,I_{n}}_{\sigma_{1},\dots,\sigma_{n}} 
= (\delta' {\cal \lambda})^{I_{1},\dots,I_{n}}_{\sigma_{1},\dots,\sigma_{n}},
\quad |I_{1}| = \dots = |I_{n}| = r-1.
\ee

Now we insert this expression into (\ref{2-3-cal}) for $l = n$ and using 
corollary \ref{delta-square} we end up with a system of the same type as above.
The recurrence establishes in the end that we have

\be
{\cal L}^{I_{1},\dots,I_{k}}_{\sigma_{1},\dots,\sigma_{k},i_{k+1},\dots,i_{n}}
 =
(\delta' {\cal \lambda})^{I_{1},\dots,I_{k}}_{\sigma_{1},\dots,\sigma_{k},
i_{k+1},\dots,i_{n}}, \quad |I_{1}| = \dots = |I_{k}| = r-1, 
\quad k = 0,\dots,n. 
\ee

It remains to substitute this relation into (\ref{L-tr}) to obtain the formula
from the statement.  
$\qed$

\begin{rem}

Because the hyper-Jacobians are traceless expressions, the $\delta$-terms
$\sum_{p=1}^{k} B_{p}\lambda$
from the expression (\ref{lagr-triv}) give a null contribution. In this way we
get the same expression as in \cite{Gr4}.
\end{rem}

We close this Subsection with a result concerning the linear dependence of the
hyper-Jacobians. 

\begin{prop}

Let
$
{\cal L}^{I_{1},\dots,I_{k}}_{\sigma_{1},\dots,\sigma_{k},i_{k+1},\dots,i_{n}},
\quad |I_{1}| = \cdots = |I_{k}| = s-1
$
be tensors depending on the variables
$(x^{i},y^{\sigma},y^{\sigma}_{j},\dots,y^{\sigma}_{j_{1},\dots,j_{s-1}})$
and having the symmetry properties (\ref{antisymmetry}). Let us suppose that
they verify the identity:

\be
\sum_{|I_{1}|=\cdots =|I_{k}|=s-1}
{\cal L}^{I_{1},\dots,I_{k}}_{\sigma_{1},\dots,\sigma_{k},i_{k+1},\dots,i_{n}}
{\cal J}^{\sigma_{1},\dots,\sigma_{k},i_{k+1},\dots,i_{n}}_{I_{1},\dots,I_{k}}
= 0
\ee

Then they are linear combinations of $\delta$-terms and conversely. In
particular, if the tensors
$
{\cal L}^{I_{1},\dots,I_{k}}_{\sigma_{1},\dots,\sigma_{k},i_{k+1},\dots,i_{n}}
$
are traceless, then the equation above has only the trivial solution.
\end{prop}

{\bf Proof:}
On denotes by $f$ the left hand side of the equation above and proves by direct
computations that:
$$
\partial^{I_{1}}_{\sigma_{1}} f = - k \sum_{|I_{2}|=\cdots =|I_{k}|=s-1}
(B_{1} {\cal L})^{I_{1},\dots,I_{k}}_{\sigma_{1},\dots,\sigma_{k},
i_{k+1},\dots,i_{n+1}}
{\cal J}^{\sigma_{2},\dots,\sigma_{k},i_{k+1},\dots,i_{n+1}}_{I_{1},\dots,
I_{k}} = 0
$$

If we continue by recurrence we obtain finally
$$
\partial^{I_{1}}_{\sigma_{1}} \cdots \partial^{I_{k}}_{\sigma_{k}} f =
(-1)^{k} k! 
(B_{1} \cdots B_{k} {\cal L})^{I_{1},\dots,I_{k}}_{\sigma_{1},\dots,\sigma_{k}}
= 0.
$$

Now one applies the basic lemma of \cite{Gr3}.
$\qed$
\newpage

\section{Two Particular Cases}

In this Section we present two particular case which do to have physical
relevance. We will obtain, essentially, known results but we think that it is
profitable to see how they follow as particularizations of the main framework
developed in this paper. There will also be some refinements of these old 
results.

\subsection{The Case $s = 2$ and $n$ arbitrary}

If we particularise in this case the main theorem of Section 3, we get:

\begin{thm}

Let $X$ be a differential manifold of dimension $N > n$ and let 
$P^{2}_{n}X$ 
be the second order Grassmann manifold associated to it. Let
$\alpha \in \Omega^{2,2}_{n+1}(X)$
be closed and verifying the Lepage condition (\ref{T2}) and the tracelessness
condition (\ref{T3}). Then $\alpha$ has the following local expression:

\begin{eqnarray}
\alpha = \sum_{k=1}^{n}
T^{i_{0}}_{\sigma_{0},\dots,\sigma_{k},i_{k+1},\dots,i_{n}} 
\omega^{\sigma_{0}}_{i_{0}} \wedge \omega^{\sigma_{1}} \wedge \cdots 
\omega^{\sigma_{k}} \wedge dx^{i_{k+1}} \wedge \cdots dx^{i_{n}} +
\nonumber \\ \sum_{k=1}^{n}
T_{\sigma_{0},\dots,\sigma_{k},i_{k+1},\dots,i_{n}} 
\omega^{\sigma_{0}} \wedge \cdots \omega^{\sigma_{k}} 
\wedge dx^{i_{k+1}} \wedge \cdots dx^{i_{n}},
\end{eqnarray}
where:

- the coefficients 
$T^{i_{0}}_{\sigma_{0},\dots,\sigma_{k},i_{k+1},\dots,i_{n}} $
and
$T_{\sigma_{0},\dots,\sigma_{k},i_{k+1},\dots,i_{n}} $
are smooth functions of the variables
$(x^{i},y^{\sigma},y^{\sigma}_{j},y^{\sigma}_{jl})$;

- they verify the symmetry properties:

\be
T^{i_{0}}_{\sigma_{0},\sigma_{P(1)},\dots,\sigma_{P(k)},
i_{Q(k+1)},\dots,i_{Q(n)}} = (-1)^{|P|+|Q|}
T^{i_{0}}_{\sigma_{0},\dots,\sigma_{k},i_{k+1},\dots,i_{n}} 
\ee

and

\be
T_{\sigma_{P(0)},\dots,\sigma_{P(k)},i_{Q(k+1)},\dots,i_{Q(n)}} =
(-1)^{|P|+|Q|} T_{\sigma_{0},\dots,\sigma_{k},i_{k+1},\dots,i_{n}};
\ee
here $P$ and $Q$ are permutations of the corresponding indices.

- the following traceless condition is valid:

\be
T^{j}_{\sigma_{0},\dots,\sigma_{k},j,i_{k+2},\dots,i_{n}} = 0.
\ee

The form $\alpha$ is globally defined by these conditions.
\end{thm}

The proof follows directly from Theorem \ref{canonical} is we take note that
the condition (\ref{particular}) is true in this particular case. Let us remark
that the previous tracelessness condition is weaker than the global condition 
$K \alpha = 0$
introduced in \cite{GP}, \cite{Gr1}.

For the sake of completeness we give below the local form of the closedness
condition 
$d\alpha = 0.$

- for $k \geq 1:$
 
\be
\partial^{lm}_{\nu} 
T^{i_{0}}_{\sigma_{0},\dots,\sigma_{k},i_{k+1},\dots,i_{n+1}} = 0,
\ee

\begin{eqnarray}
\partial^{lm}_{\sigma_{0}} 
T_{\sigma_{1},\dots,\sigma_{k},i_{k+1},\dots,i_{n+1}} +
{k+1 \over 2(n+1-k)} \sum_{p=k+1}^{n+1} (-1)^{p} \times
\nonumber \\ \left(\delta^{l}_{i_{p}}
T^{m}_{\sigma_{0},\dots,\sigma_{k},i_{k+1},\dots,\hat{i_{p}},\dots,i_{n+1}} +
\delta^{m}_{i_{p}}
T^{l}_{\sigma_{0},\dots,\sigma_{k},i_{k+1},\dots,\hat{i_{p}},\dots,i_{n+1}} 
\right) = 0, 
\end{eqnarray}

- for $k \geq 2:$

\begin{eqnarray}
{1 \over k} \left( \partial^{l}_{\sigma_{1}} 
T^{j}_{\sigma_{2},\dots,\sigma_{k},i_{k+1},\dots,i_{n+2}} -
\partial^{l}_{\sigma_{2}} 
T^{j}_{\sigma_{1},\sigma_{3},\dots,\sigma_{k},i_{k+1},\dots,i_{n+2}} \right) +
\nonumber \\
(B_{1}T)^{lj}_{\sigma_{1},\sigma_{3},\dots,\sigma_{k},i_{k+1},\dots,i_{n+2}} +
(B_{2}T)^{lj}_{\sigma_{1},\sigma_{3},\dots,\sigma_{k},i_{k+1},\dots,i_{n+2}} 
= 0
\end{eqnarray}

\begin{eqnarray}
{1 \over k} \left( \partial^{l}_{\sigma_{1}} 
T_{\sigma_{2},\dots,\sigma_{k},i_{k+1},\dots,i_{n+2}} +
\sum_{p=2}^{k} (-1)^{p-1} \partial_{\sigma_{p}} 
T^{j}_{\sigma_{1},\dots,\hat{\sigma_{p}},\dots,\sigma_{k},
i_{k+1},\dots,i_{n+2}} \right) + {1 \over n+2-k} \times  \nonumber \\
\sum_{p=k+1}^{n+2} (-1)^{p-1} d_{i_{p}} 
T^{j}_{\sigma_{1},\dots,\sigma_{k},i_{k+1},\dots,\hat{i_{p}},\dots,i_{n+2}} +
{1 \over n+2-k} \sum_{p=k+1}^{n+2} (-1)^{p-1} 
\delta^{j}_{i_{p}} 
T_{\sigma_{1},\dots,\sigma_{k},i_{k+1},\dots,\hat{i_{p}},\dots,i_{n+2}},
\end{eqnarray}
(where we assume that
$T^{i_{0}}_{i_{1},\dots,\dots,i_{n}} = 0$
according to the Lepage condition) and

\begin{eqnarray}
{1 \over k} \sum_{p=1}^{k} (-1)^{p-1} \partial_{\sigma_{p}}
T_{\sigma_{1},\dots,\hat{\sigma_{p}},\dots,\sigma_{k},i_{k+1},\dots,i_{n+2}} +
{1 \over n+2-k} \sum_{p=k+1}^{n+2} (-1)^{p-1} 
d_{i_{p}} 
T_{\sigma_{1},\dots,\sigma_{k},i_{k+1},\dots,\hat{i_{p}},\dots,i_{n+2}} 
\nonumber \\ = 0.
\end{eqnarray}

These relations are nothing but particular cases of the system (\ref{Ts}) and
(\ref{deltaT}) for $s = 2.$
As a consequence, one can prove directly that the expressions
$T_{\sigma_{0},\dots,\sigma_{k},i_{k+1},\dots,i_{n}}$
are at most linear in the second order derivatives. This can be also obtained
as a particular case of Corollary \ref{lin}. Moreover one can express
all the coefficients of the form $\alpha$ in terms of the expressions
$T_{\sigma}$;
as a consequence we have
$\alpha = 0 \Longleftrightarrow T_{\sigma} = 0.$
This is a stronger form of the relation (\ref{alpha-T}).

Finally we have the analogue of Theorem \ref{alpha-dtheta}:

\begin{thm}

In the conditions of the preceding theorem one can write $\alpha$ locally as
follows: 

\be
\alpha = d\theta
\ee
where

\be
\theta = \sum_{k=0}^{n} {1 \over k!} {n \choose k} 
\partial^{i_{1}}_{\sigma_{1}} \cdots \partial^{i_{k}}_{\sigma_{k}}
L_{i_{1},\dots,i_{n}} \omega^{\sigma_{1}} \wedge \cdots 
\omega^{\sigma_{k}} \wedge dx^{i_{k+1}} \wedge \cdots dx^{i_{n}}.
\ee

Here we have

\be
L_{i_{1},\dots,i_{n}} = \varepsilon_{i_{1},\dots,i_{n}} L
\ee
where $L$ is a smooth local function depending on the variables
$(x^{i},y^{\sigma},y^{\sigma}_{j})$.
\end{thm}

{\bf Proof:}
We have from Theorem \ref{alpha-dtheta} that the form $\theta$ has the generic
expression: 

\be
\theta = \sum_{k=0}^{n} 
L_{\sigma_{1},\dots,\sigma_{k},i_{k+1},\dots,i_{n}}
\omega^{\sigma_{1}} \wedge \cdots \omega^{\sigma_{k}} \wedge 
dx^{i_{k+1}} \wedge \cdots dx^{i_{n}}
\ee
with $L_{\dots}$ smooth local functions depending on the variables
$(x^{i},y^{\sigma},y^{\sigma}_{j},y^{\sigma}_{jl})$.

Now one has from corollary \ref{cond-theta} (or by direct commutations) the 
following relations: 

\be
\partial^{jl}_{\nu} L_{\sigma_{1},\dots,\sigma_{k},i_{k+1},\dots,i_{n}} = 0,
\ee
and for $k = 0, \dots, n$ we can express the coefficients of the form 
$\alpha$ in terms of the functions
$L_{\sigma_{1},\dots,\sigma_{k},i_{k+1},\dots,i_{n}}:$

\begin{eqnarray}
T_{\sigma_{0},\dots,\sigma_{k},i_{k+1},\dots,i_{n}} = 
{1 \over k+1} \sum_{p=0}^{k} (-1)^{p} \partial_{\sigma_{p}}
L_{\sigma_{0},\dots,\hat{\sigma_{p}},\dots,\sigma_{k},i_{k+1},\dots,i_{n}} +
\nonumber \\ {1 \over n-k} \sum_{p=k+1}^{n} (-1)^{k} d_{i_{p}}
L_{\sigma_{0},\dots,\sigma_{k},i_{k+1},\dots,\hat{i_{p}},\dots,i_{n}},
\end{eqnarray}

\begin{eqnarray}
T^{i_{0}}_{\sigma_{0},\dots,\sigma_{k},i_{k+1},\dots,i_{n}} = 
{1 \over k+1} \partial^{i_{p}}_{\sigma_{0}}
L_{\sigma_{1},\dots,\sigma_{k},i_{k+1},\dots,i_{n}} +
(B_{0}L)^{i_{0}}_{\sigma_{1},\dots,\sigma_{k},i_{k+1},\dots,i_{n}} = 
\nonumber \\
{1 \over k+1} \partial^{i_{p}}_{\sigma_{0}}
L_{\sigma_{1},\dots,\sigma_{k},i_{k+1},\dots,i_{n}} +
(-1)^{k+1} {1 \over n-k} \sum_{p=k+1}^{n} (-1)^{p-k-1} \delta^{i_{0}}_{i_{p}}
L_{\sigma_{0},\dots,\sigma_{k},i_{k+1},\dots,\hat{i_{p}},\dots,i_{n}}.
\end{eqnarray}

Now we impose the tracelessness condition from the statement of the preceding
theorem and obtain

\be
L_{\sigma_{1},\dots,\sigma_{k},i_{k+1},\dots,i_{n}} = (-1)^{k-1} 
{n-k+1 \over k^{2}} \partial ^{i_{1}}_{\sigma_{1}} 
L_{\sigma_{2},\dots,\sigma_{k},i_{1},i_{k+1},\dots,i_{n}}, 
\quad k \geq 2.
\ee

From the Lepage condition one gets easily that the condition above is true for 
$k = 1$ also. Now one obtains by induction that:

\be
L_{\sigma_{1},\dots,\sigma_{k},i_{k+1},\dots,i_{n}}  = 
{1 \over k!} {n \choose k} 
\partial^{i_{1}}_{\sigma_{1}} \cdots \partial^{i_{k}}_{\sigma_{k}}
L_{i_{1},\dots,i_{n}}
\ee
and the proof is finished.
$\qed$

We remark in the end that there are no restriction on the local Lagrangian $L$
other that the independence on the second order derivatives. 

\subsection{The Case $n = 1$ and $s$ arbitrary}

This case was studied in \cite{Kla}, \cite{Kr10}, \cite{Kr11} and \cite{Gr2}
with minor differences. As before we have:

Let $X$ be a differential manifold of dimension $N = n + 1$ and let
$P^{s}_{1}X$
be the associated Grassmann manifold of order $s$. We will denote the local
coordinates on it as follows: 
$x^{1} \mapsto t$
and
$y^{\sigma}_{\underbrace{1,\dots,1}_{k-times}} \mapsto y^{\sigma}_{k}$
i.e. we have the coordinates 
$(t,q^{\sigma}_{0},q^{\sigma}_{1},\dots,q^{\sigma}_{k})$;
it is natural to put also: 
$\omega^{\sigma}_{\underbrace{1,\dots,1}_{k-times}} \mapsto 
\omega^{\sigma}_{k}.$

\begin{thm}

In the conditions described above, let 
$\alpha \in \Omega^{s,s}_{2}(X)$
closed and verifying the tracelessness and the Lepage conditions. Then one can
write it locally as follows:

\be
\alpha = T_{\sigma} \omega^{\sigma} \wedge dt + \sum_{i+j \leq s-1} 
T^{ij}_{\sigma\nu} \omega^{\sigma}_{i} \wedge \omega^{\nu}_{j}
\label{alpha1}
\ee
with 
$T_{\sigma}, \quad T^{ij}_{\sigma\nu}$
smooth functions depending on the variables
$(t,q^{\sigma}_{0},q^{\sigma}_{1},\dots,q^{\sigma}_{k})$
verifying:

\be
T^{ij}_{\sigma\nu} = - T^{ji}_{\nu\sigma}.
\ee

In this case the form $\alpha$ is globally defined.
\end{thm}

The proof follows, as in the preceding theorem, from theorem \ref{canonical}
and the observation that the condition (\ref{particular}) is verified in this
case.  We define now the {\it total derivative operator} by:

\be
d_{t} \equiv d_{t}^{s} \equiv {\partial \over \partial t} + 
\sum_{j=0}^{s-1} q^{\sigma}_{j+1} {\partial \over q^{\sigma}_{j}}
\ee
and assume that

\be
T^{ij}_{\sigma\nu} = 0, \quad i + j \geq s.
\label{ijs}
\ee

We explicitate like before the closedness conditions for this local
description: 

\be
{\partial T_{\sigma} \over \partial q^{\nu}_{s}} + 2 T^{s-1,0}_{\nu\sigma} = 0,
\label{da1}
\ee

\be
{\partial T^{ij}_{\sigma\nu} \over \partial q^{\rho}_{s}} = 0,
\quad i,j = 0,\dots,s-1.
\label{da2}
\ee

\be
{1 \over 2} {\partial T_{\sigma} \over \partial q^{\nu}_{j}} + 
d_{t} T^{j,0}_{\nu\sigma} + T^{j-1,0}_{\nu\sigma} = 0,
\quad j = 1,\dots, s-1,
\label{da3}
\ee

\be
d_{t} T^{ij}_{\sigma\nu} + T^{i-1,j}_{\sigma\nu} + T^{i,j-1}_{\sigma\nu} = 0,
\quad i,j = 1,\dots, s-1,
\label{da4}
\ee

\be
{\partial T^{ij}_{\sigma\nu} \over \partial q^{\rho}_{k}} +
{\partial T^{jk}_{\nu\rho} \over \partial q^{\sigma}_{i}} +
{\partial T^{ki}_{\rho\sigma} \over \partial q^{\nu}_{j}} = 0,
\quad i,j,k = 0,\dots, s-1,
\label{da5}
\ee

\be
{1 \over 2} \left( {\partial T_{\sigma} \over \partial q^{\nu}_{0}} -
{\partial T_{\nu} \over \partial q^{\sigma}_{0}} \right) +
d_{t} T^{00}_{\nu\sigma} = 0.
\label{da6}
\ee

These relations follow from (\ref{Ts}) and (\ref{deltaT}) and imply in
particular that the expressions 
$T_{\sigma}$
are at most linear in the highest order derivatives
$q^{\sigma}_{s}.$

\begin{rem}

One can relax the condition (\ref{ijs}). Indeed one can accept that the form 
$\alpha$ is given by the expression (\ref{alpha1}) with the summations
restricted only to
$i,j \leq s-1.$ In that case one obtains from the closedness condition, beside
the relations (\ref{da1}) - (\ref{da6}) above, also:

\be
T^{s-1,j}_{\sigma\nu} = 0, \quad j,\dots,s-1.
\ee

Now one uses (\ref{da3}) to prove by induction that we have in fact
(\ref{ijs}) \cite{Kr11}.
\end{rem}

Let us mention two more facts. First, we have from (\ref{T_2}) in our
particular case (or directly from (\ref{da1}) - (\ref{da4}) + (\ref{da6})):

\be
T^{ij}_{\sigma\nu} = - {1 \over 2} \sum_{k=0}^{s-1-i-j} (-1)^{k+j}
{j+k \choose k} (d_{t})^{k} {\partial T^{\nu} \over \partial
q^{\sigma}_{i+j+k+1}}
\label{Tijsn}
\ee
so in particular we have 
$\alpha = 0  \Longleftrightarrow T_{\sigma} = 0.$

From this expression one can obtain the Helmholtz equations as in Corollary
\ref{helmholtz}. They are:

\be
{\partial T_{\sigma} \over \partial q^{\nu}_{j}} = 
\sum_{k=j}^{s} (-1)^{k} {k \choose j} (d_{t})^{k-j} 
{\partial T_{\nu} \over \partial q^{\sigma}_{k}}, \quad j = 0,\dots,s.
\ee

We also note that the expressions (\ref{Tijsn}) from above verify identically
the system (\ref{da1}) - (\ref{da6}). Indeed, only the equation (\ref{da5})
should be investigated because the others are used completely in the induction
process to obtain (\ref{Tijsn}). But it is not very hard to prove that
(\ref{Tijsn}) verify identically (\ref{da5}) \cite{Kr10}, \cite{Kr11}.

Finally we give the analogue of Theorem \ref{alpha-dtheta} in this case:

\begin{thm}

In the conditions of the theorem above one can write $\alpha$ locally in the form

\be
\alpha = d\theta
\ee
where

\be
\theta = \sum_{j=0}^{r-1} L^{j}_{\sigma} \omega^{\sigma}_{j}.
\ee

Here $r = [(s+1)/2]$ as before,

\be
L^{j}_{\sigma} \equiv \sum_{i=0}^{r-1-j} (-1)^{i} (d_{t})^{i} {\partial L \over
\partial q^{\sigma}_{i+j+1}}, \quad j = 0,\dots,r-1
\ee
and $L$ is a smooth function depending only on the variables:
$(t,q^{\sigma},q^{\sigma}_{j},\dots,q^{\sigma}_{r})$
which remains arbitrary for $s = 2r$ and is constrained to be at most linear in
$q^{\sigma}_{r}$ for $s = 2r-1.$
\end{thm}

The proof is elementary. We provide finally the expressions of the coefficients
of the form $\alpha$ in terms of $L$:

\be
T_{\sigma} = \sum_{j=0}^{r} (-1)^{j} (d_{t})^{j} 
{\partial L \over \partial q^{\sigma}_{j}} 
\ee
(i.e. the usual Euler-Lagrange expressions) and

\be
T^{ij}_{\sigma\nu} = {\partial L^{j}_{\nu} \over \partial q^{\sigma}_{i}} - 
{\partial L^{i}_{\sigma} \over \partial q^{\nu}_{i}}, \quad i,j = 0,\dots,s-1.
\ee

\newpage
\section{Conclusions}

We first mention that one can use the formalism developed in this paper to
analyse higher order Lagrangian systems with Noetherian symmetries, as in
\cite{Gr1}, \cite{Gr2}. Indeed, is $\phi$ is a diffeomorphisms of the manifold
$X$ then one can see that its lift 
$J^{s} \phi$
to 
$P^{s}_{n}X$
leaves invariant the subspace of forms appearing in the left hand side of
(\ref{particular}). This means that we can define a Noetherian symmetry as a
map $\phi$ such that $J^{s} \phi$ leaves the Lagrange-Souriau class invariant.
It is to be expected that the computations will be much more difficult than in
the two particular cases from the last Section.

Next, we mention that it is not clear if in the general case studied here, the 
only restriction on the Euler-Lagrange expressions are given by the generalized
Helmholtz equations, but it is reasonable to conjecture that this is true.

Last, we remark that the formalism above could be generalised, in principle, to
the case when the Euler-Lagrange expressions are not restricted by the
condition of linearity in the highest order derivatives, trying for instance to
relax the condition (\ref{T1})  i.e to factorize $\alpha$ to a smaller 
subspace.

%{\bf Acknowledgments:} 


\newpage
\begin{thebibliography}{99}

\bibitem{AM}
H. F. Ahner, A. E. Moose,
``{\it Covariant Inverse Problem of the Calculus of Variation}",
Journ. Math. Phys. {\bf 18} (1977) 1367-1373
 
\bibitem{An}
I. M. Anderson,
``{\it The Variational Bicomplex}",
Utah State Univ. preprint, 1989, (Academic Press, Boston, to appear)

\bibitem{AD}
I. M. Anderson, T. Duchamp,
``{\it On the Existence of Global Variational Principles}",
American Journ. Math. {\bf 102} (1980) 781-868

\bibitem{AH}
R. W. Atherton, G. M. Homsy,
``{\it On the Existence and Formulation of Variational Principles for
Nonlinear Differential Equations}",
Studies in Appl. Math. {\bf LIV} (1975) 31-60

\bibitem{Be1} 
D. E. Betounes, 
``{\it Extensions of the Classical Cartan Form}", 
Phys. Rev. {\bf D 29} (1984) 599-606

\bibitem{Be2} 
D. E. Betounes, 
``{\it Differential Geometric Aspects of the Cartan Form: Symmetry Theory}", 
J. Math. Phys. {\bf 28} (1987) 2347-2353

\bibitem{BCO}
J. M. Ball, J. C. Currie, P. J. Olver,
``{\it Null Lagrangians, Weak Continuity, and Variational Problems
of Arbitrary Order}",
Journ. Functional Anal. {\bf 41} (1981) 135-174

\bibitem{Ca} 
E. Cartan, 
``{\it Le\c cons sur les Invariants Integraux}", 
Hermann, 1922.

\bibitem{Ga} 
P. L. Garcia, 
``{\it The Poincar\'e-Cartan Invariant in the Calculus of Variations}", 
Symp. Math. {\bf XIV} (1974) 219-246

\bibitem{Go} 
M. J. Gotay, 
``{\it A Multisymplectic Framework for Classical Field Theory and the 
Calculus of Variations. I. Covariant Hamiltonian  Formalism}", 
in ``Mechanics, Analysis and Geometry: 200 Years after Lagrange", 
M. Francaviglia and D. D. Holms, eds., North-Holland, Amsterdam, 1990, 
pp. 203-235

\bibitem{Gr1}
D. R. Grigore,
``{\it A Generalized Lagrangian Formalism in Particle Mechanics and Classical
Field Theory}",
Fortschr. der Phys. {\bf 41} (1993) 569-617

\bibitem{Gr2}
D. R. Grigore,
``{\it Higher-Order Lagrangian Theories and Noetherian Symmetries}",
Romanian Journ. Phys. {\bf 39} (1994) 11-35

\bibitem{Gr3}
D. R. Grigore,
``{\it Variational Sequence on Finite Jet Bundle Extensions and the Lagrangian 
Formalism}",
dg-ga/9702016, submitted for publication

\bibitem{Gr4} 
D. R. Grigore,
``{\it Variationally Trivial Lagrangians and Locally Variational Differential 
Equations of Arbitrary Order}", 
submitted for publication, 

\bibitem{Gr5} 
D. R. Grigore,
``{\it Higher-Order Lagrangian Formalism on Grassmann Manifolds}", 
dg-ga/9709005, submitted for publication, 

\bibitem{GK} 
D. R. Grigore, D. Krupka,
``{\it Invariants of Velocities and Higher Order Grassmann Bundles}",
Journ. Geom. Phys. {\bf 24} (1997) 244-266, dg-ga/9708013

\bibitem{GM}
A. Galindo, L. Mart\'\i nez Alonso,
``{\it Kernels and Ranges in the Variational Formalism}",
Lett. Math. Phys. {\bf 2} (1978) 385-390

\bibitem{GP} 
D. R. Grigore, O. T. Popp, 
``{\it On the Lagrange-Souriau Form in Classical Field Theory}", 
to appear in Mathematica Bohemica

\bibitem{GS} 
H.Goldschmits, S.Sternberg: 
``{\it The Hamilton-Cartan Formalism in the Calculus of Variations}", 
Ann. Inst. Fourier {\bf 23} (1973) 203-267

\bibitem{Horv} 
P. Horv\`athy,
``{\it Variational Formalism for Spinning Particles}",
Journ. Math. Phys. {\bf 20} (1979) 49-52

\bibitem{Ki} 
I. Kijowski, 
``{\it A Finite-Dimensional Canonical Formalism in Classical Field Theory}", 
Comm. Math. Phys. {\bf 30} (1973) 99-128

\bibitem{Kla}
L. Klapka,
{\it Euler-Lagrange Expressions and Closed Two-Forms in Higher Order
Mechanics}, 
in {\it Geometrical Methods in Physics}, Conf. on Differential Geometry and
Applications, Czechoslovakia, 1983, (Univ. Brno, D. Krupka ed.) 

\bibitem{Kl} 
J. Klein, 
``{\it Espaces Variationels et M\'ecaniqu}e", 
Ann. Inst. Fourier {\bf 12} (1962) 1-124

\bibitem{Kr1} 
D. Krupka, 
``{\it A Map Associated to the Lepagean Forms of the Calculus 
of Variations in Fibered Manifolds}", 
Czech. Math. Journ. {\bf 27} (1977) 114-118

\bibitem{Kr2}
D. Krupka,
``{\it Lepagean Forms in Higher Order Variational Theory}",
in Proceedings of the IUTAM-ISIMM Symposium on ``Modern Developments in
Analytical Mechanics", Turin, 1982, Atti della Academia delle Scienze di
Torino, Suppl. al Vol. {\bf 117} (1983) 198-238

\bibitem{Kr3}
D. Krupka,
``{\it Geometry of Lagrangian Structures.3}",
in Proceedings of the 14th Winter School on Abstract Analysis, SRNI, 1986,
Suppl. ai Rendiconti del Circolo Matematico di Palermo, Serie II, no.14, 1987,
pp. 187-224

\bibitem{Kr4}
D. Krupka,
``{\it Variational Sequence on Finite Order Jet Spaces}",
in Proceedings of the Conference ``Differential Geometry and its Applications",
August, 1989, World Scientific, Singapore, 1990, pp. 236-254

\bibitem{Kr5}
D. Krupka,
``{\it Topics in the Calculus of Variation: Finite Order Variational
Sequences}", 
in ``Differential Geometry and its Applications", proceedings Conf. Opava,
1992, Open Univ. Press, pp. 437-495

\bibitem{Kr9}
D. Krupka,
``{\it The Geometry of Lagrange Structure}",
Preprint Series in Global Analysis, GA 7/97, Dept. Of Math., Opava Univ.,
Czech Rep.

\bibitem{Kr10}
Olga Krupkov\'a,
{\it Lepagean 2-forms in Higher Order Hamiltonian Mechanics. II Inverse
Problem}, 
Arch. Math. (Brno) {\bf 22} (1986) 97-129

\bibitem{Kr11}
Olga Krupkov\'a,
{\it Variational Analysis on Fibered Manifolds over One-Dimensional Bases},
Ph. D. Thesis, Opava Univ., 1992

\bibitem{Ol1}
P. J. Olver,
``{\it Hyperjacobians, Determinant Ideals and the Weak Solutions to Variational
Problems}",
Proc. Roy. Soc. Edinburgh {\bf 95A} (1983) 317-340

\bibitem{Ol2}
P. J. Olver,
``{\it Applications of Lie Groups to Differential Equations},
Springer, 1986

\bibitem{Po} 
H. Poincar\'e, 
``{\it Le\c cons sur les M\'ethodes Nouvelles de la  M\'ecanique C\'eleste}", 
Gauthier-Villars, Paris, 1892.

\bibitem{Ru1} 
H. Rund, 
``{\it A Cartan Form for the Field Theory of Charath\'eodory in 
the Calculus of Variations of Multiple Integrals}", in
Lect. Notes in Pure and Appl. Math. {\bf 100} (1985) 455-469

\bibitem{Ru2} 
H. Rund,
``{\it Integral Formulae Associated with the Euler-Lagrange Operator
of Multiple Integral Problems in the Calculus of Variation}",
\AE quationes Math. {\bf 11} (1974) 212-229

\bibitem{Sa} 
D. J. Saunders, 
``{\it An Alternative Approach to the Cartan Form in the Lagrangian Field 
Theories}", 
J. Phys.. {\bf A 20} (1987) 339-349

\bibitem{Sou} 
J. M. Souriau, 
``{\it Structure des Systemes Dynamique}", 
Dunod, Paris, 1970.

\end{thebibliography}
\end{document}